\documentclass[11pt]{article}
\usepackage{amssymb,amscd,amsthm,tikz-cd,mathrsfs,textcomp,mathtools,setspace,titlesec,url,hyperref}

\titleformat{\subsection}[runin]
  {\normalfont\bfseries}{\thesubsection}{1em}{}

\RequirePackage[left=1in,right=1in,top=1in,bottom=1in]{geometry}

\newcommand{\Z}{\mathbb{Z}}

\newcommand{\A}{\mathbb{A}}
\newcommand{\G}{\mathbb{G}}

\DeclareMathOperator{\AJ}{AJ}

\DeclareMathOperator{\Gr}{Gr}
\DeclareMathOperator{\id}{id}

\DeclareMathOperator{\dR}{dR}

\DeclareMathOperator{\Hom}{Hom}
\DeclareMathOperator{\Map}{Map}
\DeclareMathOperator{\Ext}{Ext}
\DeclareMathOperator{\ext}{ext}
\DeclareMathOperator{\red}{red}
\DeclareMathOperator{\aff}{aff}
\DeclareMathOperator{\can}{can}

\DeclareMathOperator{\add}{add}

\DeclareMathOperator{\hor}{hor}
\DeclareMathOperator{\Rat}{Rat}
\DeclareMathOperator{\Ran}{Ran}
\DeclareMathOperator{\Spec}{Spec}

\DeclareMathOperator{\univ}{univ}

\DeclareMathOperator{\colim}{colim}

\DeclareMathOperator{\uHom}{\underline{Hom}}
\DeclareMathOperator{\uMap}{\underline{Map}}
\DeclareMathOperator{\Pic}{Pic}
\DeclareMathOperator{\Proj}{Proj}

\DeclareMathOperator{\Div}{Div}
\DeclareMathOperator{\ddiv}{div}
\DeclareMathOperator{\disj}{disj}
\DeclareMathOperator{\uExt}{\underline{Ext}}

\numberwithin{equation}{subsection}

\newtheorem{theorem}{Theorem}[subsection]
\newtheorem{corollary}{Corollary}[theorem]
\newtheorem{lemma}[theorem]{Lemma}
\newtheorem{proposition}[theorem]{Proposition}
\newtheorem{corollary'}[theorem]{Corollary}

\theoremstyle{definition}
\newtheorem{definition}[theorem]{Definition}

\title{Geometric class field theory and Cartier duality}
\author{Justin Campbell and Andreas Hayash}

\begin{document}

\maketitle

\begin{abstract}
We formulate and prove a generalized Albanese property for families of maps from a smooth curve over an arbitrary field into a commutative group stack. Our proof, which is mostly self-contained, employs local-to-global techniques and some new Ext-vanishing results to reduce to the local Cartier self-duality theorem of Contou-Carr\'{e}re. As a corollary, we reprove local and global geometric class field theory with arbitrary ramification.
\end{abstract}

\section{Introduction}

\subsection{History of the problem.} Modern geometric class field theory has its origins in the work of Lang and Rosenlicht in the 1950s, summarized in the book by Serre \cite{S1}. They reproved Artin reciprocity for function fields using methods from algebraic geometry, emphasizing the role of the Albanese property enjoyed by the Jacobian variety of an algebraic curve.

Later, Deligne and Serre reinterpreted geometric class field theory in terms of \'{e}tale local systems, so that it could be formulated for curves over any field. The Artin reciprocity isomorphism for function fields over a finite ground field can be recovered from their formulation by passing through the \emph{faisceaux-fonctions} correspondence. Deligne gave a beautiful geometric proof of the global unramified case, which is to say for local systems on a smooth projective curve (see \cite{La} for an account). After some modification his argument applies to the tamely ramified case (\cite{Te},\cite{To}), and recently has been extended to allow arbitrary ramification (\cite{Gu},\cite{Ta}).

Meanwhile, local geometric class field theory is not particularly well-documented in the literature. Serre \cite{S2} gave an approach to local class field theory for a complete non-archimedean field with algebraically closed residue field of positive characteristic, which was generalized by Hazewinkel \cite{H} to allow any perfect residue field of positive characteristic, and later Suzuki and Yoshida \cite{SY} extended this to an arbitrary perfect residue field (see also \cite{V}). This approach has a geometric flavor, with the group of local units in some sense playing the role of local Jacobian. However, the arguments are arithmetic in nature, relying on Galois cohomology and related techniques.

\subsection{The generalized Albanese property.} In this paper we formulate and prove a generalized Albanese property in the spirit of Lang and Rosenlicht, valid for maps from a smooth curve over any field into a rather general kind of commutative group stack, which is required only to satisfy a reflexivity property with respect to Cartier duality (Corollary \ref{maincor}). As a consequence, we obtain a strong form of geometric class field theory \emph{\'{a} l\`{a}} Deligne, allowing local systems with arbitrary ramification (Corollary \ref{gcft}). The Albanese property applies to families of maps over an arbitrary base; in fact, this is crucial for our proof (we work with a curve defined over the ground field, but our methods should extend to the case of a relative curve without too much difficulty).

Cartier duality plays a key role in our approach. Namely, when the target stack is taken to be the classifying stack of the multiplicative group, the generalized Albanese property amounts to a perfect pairing between the Picard stack of the open curve (suitably interpreted) and the stack classifying line bundles on its projective closure with full level structure at the punctures (Theorem \ref{mainthm}). Our proof of the Albanese property in general proceeds by reduction to this case.

Although many special cases were previously known, and the subject is an old one, our main result appears to be new. One possible explanation for this is that the generalized Albanese property for open curves is subtle to formulate in families. The obstruction can be seen in some examples constructed by Bass in 1962 \cite{Bas}, which show that a family of line bundles on an open curve need not extend to its projective closure, even locally on the base. We overcome this difficulty, essentially by systematically excluding such non-extendable families of line bundles. In addition, we show that our Albanese property specializes to the na\"{i}ve one when the target is a commutative affine group scheme or the classifying stack of a finite abelian group.

\subsection{Methods.} We begin by establishing a local Cartier self-duality theorem (Theorem \ref{locramthm}) in the style of Contou-Carr\`{e}re \cite{CC1} using the theory of Tate vector bundles and the Beilinson-Kapranov determinant gerbe (see \cite{Dr},\cite{K}, and \cite{BBE} for background on these notions; this approach is compared to Contou-Carr\`{e}re's in Appendix \ref{appendix}). Our result strengthens Contou-Carr\`{e}re's by also proving that any family of multiplicative line bundles on the group of local units is trivial locally on the base. In fact, we prove a similar vanishing statement for an arbitrary commutative affine group scheme, which surprisingly seems to be new even for the additive group, although many other cases were known previously (see \cite{Br}).

Our proof of the local self-duality theorem is carried out factorizably, meaning we consider several points of the curve simultaneously and allow them to move and collide. This allows us to use local-to-global methods when proving the global Cartier self-duality theorem. The other main input for the global theorem is the fact that any family of multiplicative line bundles on the group of global invertible rational functions is trivial locally on the base (Lemma \ref{vanish}; see \cite{BD} Proposition 4.9.1 for a closely related statement). The proof of this global vanishing result in \emph{loc. cit.} contains a gap, which we address following a suggestion made by Beilinson to the first-named author.

As previously mentioned, we also show that our formulation of the generalized Albanese property implies geometric class field theory with arbitrary ramification (Theorem \ref{finext}). The content of this implication is that any family of torsion line bundles on the open curve extends to its projective completion locally on the base. As we explain in Section \ref{cftsec}, this follows from recent results of  Bouthier and \v{C}esnavi\v{c}ius \cite{BC}.

\subsection{Relation to previous work.} In addition to the above-cited works, the formulation of our main theorem was inspired by \cite{De}, where a version of the generalized Albanese property is proved in the unramified setting (i.e. for a projective curve). Their hypotheses on the target commutative group stack are rather complicated and seem to be less general than ours, although we are uncertain of the precise relation.

Contou-Carr\`{e}re proved a version of the Albanese property for an open curve in \cite{CC2}, allowing families of maps into a smooth commutative group scheme. His approach, although different from ours, also emphasizes the role of Cartier duality.

Our approach has also been influenced by Mirkovi\'{c}’s conjectural inner homology theory in algebraic geometry (see \cite{M}). Specifically, he has proposed dual homology and cohomology theories on the category of (ind)-schemes with values in commutative group stacks which refine the usual singular theories. Our formulation of geometric class field theory, and in particular the generalized Albanese property, may be viewed as saying that the value of his cohomology theory on a curve is the Picard stack. In particular, the ramified duality statement may be construed as a version of Poincare duality between compactly supported cohomology and homology. One may adopt a similar perspective locally by thinking of the affine Grassmannian as the compactly supported cohomology of a disc.

Let us also comment on the relationship with the geometric Langlands program. Although geometric class field theory can be viewed as a special case of the geometric Langlands correspondence, proposed proofs of the latter (see \cite{Ga2} for a proof sketch in the de Rham setting) take class field theory as an input via the theory of geometric Eisenstein series. Thus the geometric Langlands program does not so much subsume geometric class field theory as rely on it.

We became aware of \cite{Gu} and \cite{Ta} during the preparation of this manuscript. They both prove versions of our Corollary \ref{gcft}: Takeuchi proves the isomorphism at the level of geometric points, and Guignard in families (even for a relative curve). Their methods are quite different from ours.

\subsection*{Acknowledgements} The authors thank Dennis Gaitsgory, Ivan Mirkovi\'{c}, Sam Raskin, Kyle Ferendo, and Ka Fai Wong for many helpful conversations without which this paper would not have been possible. Part of this work was carried out at the Nisyros Conference on Automorphic Representations in July 2019, and the authors also express their gratitude to the organizers.

\section{Statement of results}

\subsection{} The letters $S$, $T$, etc. will always denote affine schemes over $k$. We work with stacks in the fppf topology unless otherwise specified. All Picard groupoids and commutative group stacks are assumed to be strictly commutative, i.e. for any point $g$ the commutativity constraint $g \otimes g \tilde{\to} g \otimes g$ is required to be the identity map.

Recall that the 2-category of (strictly) commutative group stacks is equivalent to the full subcategory of the derived category of sheaves of abelian groups consisting of complexes concentrated in degrees $-1$ and $0$. Accordingly, we can speak of exact triangles of commutative group stacks, etc. See \cite{Br} for more details.

Given prestacks $Y$ and $Z$, we write $\uMap(Y,Z)$ for the internal mapping space, defined by \[ \Map(S,\uMap(Y,Z)) := \Map(Y \times S,Z). \] Similarly, if $G$ and $H$ are commutative group stacks then we denote by $\uHom(G,H)$ their internal Hom, given on the level of points by \[ \Map(S,\uHom(G,H)) := \Hom_S(G \times S,H \times S). \] In particular, we write \[ G^* := \uHom(G,\G_m) \] and \[ G^{\vee} := \uHom(G,B\G_m) \] for the Cartier $0$-dual and $1$-dual, respectively, of $G$. We say that $G$ is $0$\emph{-reflexive}, respectively $1$\emph{-reflexive}, if the canonical map $G \to G^{**}$, respectively $G \to G^{\vee \vee}$, is an isomorphism.

\subsection{} \label{extdefsec}
Let $\overline{X}$ be a smooth, projective, and geometrically connected curve over a field $k$. Fix a nonempty open subscheme $X \subset \overline{X}$. Write $\Pic(X,\partial X)$ for the commutative group stack which parameterizes line bundles on $X$ equipped with a trivialization along the ``boundary," i.e. the punctured formal neighborhood of the closed points in $\overline{X} \setminus X$.

The stack $\Pic(X,\partial X)$ receives the Abel-Jacobi map
\begin{align*}
\AJ_X : X &\longrightarrow \Pic(X,\partial X) \\
x &\mapsto (\mathcal{O}_X(x),1).
\end{align*}

We define $\Pic(X)$ to be the stack parameterizing line bundles on $X$ which extend to $\overline{X}$ locally on the base: we call such line bundles \emph{extendable}. Equivalently $\Pic(X)$ is the stack-theoretic image of the forgetful map \[ \Pic(X,\partial X) \longrightarrow \uMap(X,B\G_m). \] Of course, if $X$ is projective then $\Pic(X)$ is the usual Picard stack. Note that if $S$ is regular, then all line bundles on $X \times S$ are extendable.

Let $G$ be a $1$-reflexive commutative group stack. A map $X \times S \to G$ will be called $B\G_m$\emph{-extendable} if, for any map $T \to S$ and any homomorphism $G \times T \to B\G_m \times T$ over $T$, the resulting line bundle \[ X \times T \longrightarrow G \times T \longrightarrow B\G_m \] is extendable. We write \[ \uMap(X,G)^{\ext} \subset \uMap(X,G) \] for the substack consisting of $B\G_m$-extendable maps. Observe that we have a tautological identification \[ \uMap(X,G)^{\ext} = \uHom(G^{\vee},\Pic(X)). \] In particular $\uMap(X,B\G_m)^{\ext} = \Pic(X)$.

We can now state our main theorems.

\begin{theorem}
\label{mainthm}
Restriction along $\AJ_X$ induces an isomorphism \[ \Pic(X,\partial X)^{\vee} \tilde{\longrightarrow} \Pic(X). \] Moreover, the commutative group stacks $\Pic(X,\partial X)$ and $\Pic(X)$ are $1$-reflexive.
\end{theorem}

See Section \ref{mainthmproof} for the proof. Note that the theorem says in particular that we have a perfect pairing \[ \langle \, , \, ]_X : \Pic(X,\partial X) \times \Pic(X) \longrightarrow B\G_m. \]

The generalized Albanese property, stated as follows, follows from Theorem \ref{mainthm} by a formal argument.

\begin{corollary}
\label{maincor}
For any $1$-reflexive commutative group stack $G$, restriction along $\AJ_X$ induces an isomorphism \[ \uHom(\Pic(X,\partial X),G) \tilde{\longrightarrow} \uMap(X,G)^{\ext}. \]
\end{corollary}

This corollary is also proved in Section \ref{mainthmproof}. Note that the case $G = B\G_m$ recovers the theorem.

\subsection{} Our proof of Theorem \ref{mainthm} is by local-to-global methods. Fixing a closed point $x$ in $X$, we consider the formal disk $D_x$ and the punctured disk $\mathring{D}_x$ centered at $x$.

For any prestack $Y$, we will write $\aff(Y) := \Spec \Gamma(Y,\mathcal{O}_Y)$. Let \[ \mathring{D}_{x,S} := \aff(D_x \times S) \setminus (\{ x \} \times S). \] We abuse notation by writing $\uMap(\mathring{D}_x,Y)$ for the modified mapping prestack defined by \[ \Map(S,\uMap(\mathring{D}_x,Y)) := \Map(\mathring{D}_{x,S},Y). \]

The role of Jacobian of $\mathring{D}_x$ is played by the group \[ \mathfrak{L}_x\G_m := \uMap(\mathring{D}_x,\G_m) \] of formal loops into $\G_m$ centered at $x$. We write \[ \Pic(\mathring{D}_x) := B\mathfrak{L}_x\G_m, \] so analogously to the global situation we have an inclusion \[ \Pic(\mathring{D}_x) \longrightarrow \uMap(\mathring{D}_x,B\G_m) \] whose image consists of families satisfying an extendability condition.

The group stacks $\mathfrak{L}_x\G_m$ and $\Pic(\mathring{D}_x)$ generalize in a straightforward fashion to the case where $x$ is replaced by a finite subset of $X$, and continue to make sense in families as the points move and collide. Moreover, these spaces admit natural factorization structures which reflect their local nature. In Section \ref{localsec} we recall how this structure is encoded using the Ran space of $X$.

The key local input to our theorem is the pairing \[ \langle \, , \, ]_{\mathring{D}_x} : \Pic(\mathring{D}_x) \times \mathfrak{L}_x\G_m \longrightarrow B\G_m, \] obtained by delooping the $2$-cocyle
\begin{equation}
\label{ccpairing}
\mathfrak{L}_x\G_m \times \mathfrak{L}_x\G_m \longrightarrow \G_m
\end{equation} attached to the Tate extension of $\mathfrak{L}_x\G_m$. In Section \ref{localsec} we review the construction of this pairing using determinant gerbes of Tate vector spaces.

Note that if $I$ denotes the set of closed points in $\overline{X} \setminus X$, then we have an exact triangle \[ \Pic(X,\partial X) \longrightarrow \Pic(X) \longrightarrow \prod_{x \in I} \Pic(\mathring{D}_x). \]

\begin{theorem}
\label{locramthm}
The pairing $\langle \, , \, ]_{\mathring{D}_x}$ is perfect, and is moreover compatible with the factorization structures on $\mathfrak{L}_x\G_m$ and $\Pic(\mathring{D}_x)$. The isomorphism induced by $\langle \, , \, ]_{\mathring{D}_x}$ fits into an isomorphism of exact triangles
\begin{equation}
\label{locglobdiag}
\begin{tikzcd}
\prod_{x \in I} \mathfrak{L}_x\G_m \arrow{r} \arrow{d} & \Pic(X,\partial X) \arrow{r} \arrow{d} & \Pic(X) \arrow{d} \\
\prod_{x \in I} \Pic(\mathring{D}_x)^{\vee} \arrow{r} & \Pic(X)^{\vee} \arrow{r} & \Pic(X,\partial X)^{\vee},
\end{tikzcd}
\end{equation}
where the middle and right vertical arrows are induced by $\langle \, , \, ]_X$.
\end{theorem}

See Theorem \ref{locnondeg} below for the nondegeneracy of the factorizable pairing, and Section \ref{locramthm} for the proof of the local-global compatibility.

Note that Theorem \ref{locramthm} implies that (\ref{ccpairing}) is perfect, and also that $\underline{\Ext}^1(\mathfrak{L}_x\G_m,\G_m) = 1$.

Similarly to the case of an affine curve, to any $1$-reflexive commutative group stack $G$ we attach a stack $\uMap(\mathring{D}_x,G)^{\ext}$ of $B\G_m$-extendable maps, defined so that \[ \uMap(\mathring{D}_x,G)^{\ext} = \uHom(G^{\vee},\Pic(\mathring{D}_x)). \]

\begin{corollary}
\label{locramalb}
There is a unique natural isomorphism \[ \uHom(\mathfrak{L}_x\G_m,G) \tilde{\longrightarrow} \uMap(\mathring{D}_x,G)^{\ext}, \] where both sides are viewed as functors of $1$-reflexive group stacks $G$, which agrees with the isomorphism induced by $\langle \, , \, ]_{\mathring{D}_x}$ in the case $G = B\G_m$. For fixed $G$, this isomorphism extends to an isomorphism of factorization spaces when $x$ is allowed to vary.
\end{corollary}

See Corollary \ref{locramalb2} for a more precise formulation and its proof.

Our characterization of the isomorphism in Corollary \ref{locramalb} is more complicated than that of the analogous global isomorphism in Corollary \ref{maincor}, since the former is not given by restriction along an Abel-Jacobi map \[ \mathring{D}_x \longrightarrow \mathfrak{L}_x\G_m. \] However, such a map can be constructed after modifying $\mathfrak{L}_x\G_m$ suitably. This was the approach of Contou-Carr\`{e}re, which we compare with ours in the appendix. The construction is somewhat awkward: the authors expect that a conceptual explanation will probably involve some kind of rigid-analytic or adic geometry.

\subsection{} In order to apply our results to geometric class field theory, we take $G = BA$ where $A$ is a constant finite abelian group. It turns out the the $B\G_m$-extendability condition is superfluous in this case.

\begin{theorem}
\label{finext}
The inclusions \[ \uMap(X,BA)^{\ext} \longrightarrow \uMap(X,BA) \] and \[ \uMap(\mathring{D}_x,BA)^{\ext} \longrightarrow \uMap(\mathring{D}_x,BA) \] are isomorphisms.
\end{theorem}

See Section \ref{finextproof} for the proof. Combining the theorem with Corollary \ref{maincor}, we immediately obtain the following version of global geometric class field theory with arbitrary ramification.

\begin{corollary}
\label{gcft}
Restriction along $\AJ_X$ induces an isomorphism \[ \uHom(\Pic(X,\partial X),BA) \tilde{\longrightarrow} \uMap(X,BA). \]
\end{corollary}

Similarly, Theorem \ref{finext} combined with Corollary \ref{locramalb} yields local geometric class field theory.

\begin{corollary}
\label{lcft}
When $G = BA$, the isomorphism of Corollary \ref{locramalb} specializes to an isomorphism \[ \uHom(\mathfrak{L}_x\G_m,BA) \tilde{\longrightarrow} \uMap(\mathring{D}_x,BA), \] natural in $A$.
\end{corollary}

Passing to $k$-points, we in particular obtain isomorphisms \[ \Hom(\Pic(X,\partial X),BA) \tilde{\longrightarrow} \Map(X,BA) \] and \[ \Hom(\mathfrak{L}_x\G_m,BA) \tilde{\longrightarrow} \Map(\mathring{D}_x,BA), \] natural in $A$. In the case where $k$ is finite, standard arguments using the trace of Frobenius construction then imply the classical Artin reciprocity isomorphisms for the function field of $X$ and its completions.

\section{Tate vector bundles}

\subsection{} Below we recall the basic theory of Tate vector bundles and their determinant gerbes, developed in \cite {BBE} and \cite{Dr}. Writing $S = \Spec A$, we will consider topological $A$-modules which are complete, Hausdorff, and admit a base of open neighborhoods of $0$ consisting of $A$-submodules. Given a homomorphism $A \to B$, restriction of scalars preserves such topological modules, and moreover admits a left adjoint given by the completed extension of scalars $M \mapsto B \widehat{\otimes}_A M$.

A \emph{discrete Tate vector bundle} on $S$ is a projective $A$-module. A topological $A$-module is called a \emph{linearly compact Tate vector bundle} on $S$ if it is isomorphic to the topological $A$-dual of a projective $A$-module. An \emph{elementary Tate vector bundle} on $S$ is a topological $A$-module isomorphic to one of the form $P \oplus Q$ where $P$ and $Q$ are discrete and linearly compact, respectively. A \emph{Tate vector bundle} on $S$ is a topological $A$-module which is a direct summand of an elementary Tate vector bundle. It is shown in \cite{Dr} that any Tate vector bundle is elementary locally in the Nisnevich topology (in particular, locally in the \'{e}tale topology).

A \emph{lattice} in a Tate vector bundle $V$ is an open and linearly compact submodule $L \subset V$ such that $V/L$ is projective (in \emph{loc. cit.} it is called a coprojective lattice). Note that a Tate vector bundle is elementary if and only it admits a lattice.

\subsection{} A key role is played by the following construction from \cite{BBE}. Consider the category $\Proj(S)^{\infty}$ whose objects are projective $A$-modules and whose morphisms are given by the formula \[ \Hom_{\Proj(S)^{\infty}}(P,Q) := \Hom_A(P,Q)/\Hom_A^f(P,Q), \] where $\Hom_A^f(P,Q)$ consists of $A$-linear maps $P \to Q$ whose image is contained in a finitely generated submodule. Denoting by $\Proj(S)$ the category of projective $A$-modules, we therefore have a canonical functor $\Proj(S) \to \Proj(S)^{\infty}$.

A morphism $\varphi : P \to Q$ in $\Proj(S)$ is called \emph{Fredholm} if its image $\varphi^{\infty}$ in $\Proj(S)^{\infty}$ is an isomorphism. In \emph{loc. cit.} Proposition 2.7, the authors construct a determinant line bundle $\det(\varphi^{\infty})$ on $S$ (there denoted by $\det(Q,P,\varphi^{\infty})$) attached functorially to any Fredholm operator $\varphi$. This line bundle depends only on $\varphi^{\infty}$, its formation is compatible with base change, and if $P,Q$ have finite rank then we have $\det(\varphi^{\infty}) = (\det P)^{-1} \otimes \det Q$.

This allows us to construct the \emph{relative determinant line} of two lattices $L_1,L_2 \subset V$. Write $\pi : V \to V/L_2$ for the projection and choose a splitting $\sigma : V/L_1 \to V$ of the other projection. Then $(\pi \circ \sigma)^{\infty}$ does not depend on the choice of $\sigma$, and we put \[ \det(L_1,L_2) := \det((\pi \circ \sigma)^{\infty}), \] a line bundle on $S$.

If $L_1 \supset L_2$ then $L_1/L_2$ is projective of finite rank and $\det(L_1,L_2) = \det(L_1/L_2)$. For any three lattices $L_1,L_2,L_3 \subset V$ we have a canonical isomorphism \[ \det(L_1,L_2) \otimes \det(L_2,L_3) \tilde{\longrightarrow} \det(L_1,L_3), \] and these isomorphisms satisfy the obvious cocycle condition.

\subsection{} Following \cite{Dr}, to any Tate vector bundle $V$ we attach its determinant gerbe $D(V)$, a $\G_m$-gerbe (a.k.a. $B\G_m$-torsor) on $S$. A section of $D(V)$ is a (weak) determinant theory $d$ on $V$, which assigns to every lattice $L \subset V$ a line bundle $d(L)$ equipped with isomorphisms \[ d(L_1) \otimes \det(L_1,L_2) \tilde{\longrightarrow} d(L_2), \] compatible with the transitivity isomorphisms for the relative determinant. Formation of the determinant gerbe is compatible with base change.

Any lattice $L \subset V$ determines a determinant theory $L' \mapsto \det(L,L')$ on $V$, and hence a section of $D(V)$. In particular the gerbe $D(V)$ is trivial Nisnevich locally on $S$.

The construction $V \mapsto D(V)$ is functorial in isomorphisms. In particular, if $\varphi$ is a continuous automorphism of the Tate vector bundle $V$, it determines an automorphism $D(\varphi;V)$ of $D(V)$, which is just a line bundle on $S$. If $L \subset V$ is a lattice then we have \[ \textstyle D(\varphi;V) = \det(\varphi(L),L). \] Moreover, we have canonical isomorphisms \[ \textstyle D(\varphi_1 \circ \varphi_2;V) \tilde{\longrightarrow} D(\varphi_1;V) \otimes D(\varphi_2;V) \] satisfying the obvious compatibilities.

\subsection{} If $V$ is a Tate vector bundle on $S = \Spec A$, we call a submodule $M \subset V$ a \emph{colattice} if it is complementary to a lattice. In particular $M$ is discrete and projective, and $V/M$ is linearly compact.

We would like to attach to a pair of colattices $M_1,M_2 \subset V$ their relative determinant line bundle $\det(M_1,M_2)$ on $S$, in much the same way as for lattices. In particular we desire that if $M_1 \supset M_2$ then we have $\det(M_1,M_2) = \det(M_1/M_2)$, and that we have transitivity isomorphisms satisfying a cocycle condition.

Denote by $\iota : M_1 \to V$ the inclusion and choose a retract $\rho : V \to M_2$ of the other inclusion. Then $(\rho \circ \iota)^{\infty}$ does not depend on the choice of $\rho$, and we define \[ \det(M_1,M_2) := \det((\rho \circ \iota)^{\infty}). \]

Finally, we also introduce the relative determinant line of a lattice $L \subset V$ and a colattice $M \subset V$. Denote by $\iota : M \to V$ the inclusion and $\pi : V \to V/L$ the projection. We define \[ \det(M,L) := \det((\pi \circ \iota)^{\infty}) \] and \[ \det(L,M) := \det(M,L)^{-1}. \] Note that we have transitivity isomorphisms \[ \det(M,L) \otimes \det(L,L') \tilde{\longrightarrow} \det(M,L') \] and \[ \det(M',M) \otimes \det(M,L) \tilde{\longrightarrow} \det(M',L) \] satisfying the obvious compatibility conditions. In particular a choice of colattice determines a determinant theory on $V$, i.e. trivializes the determinant gerbe $D(V)$.

\section{The local pairing}

\label{localsec}

\subsection{} Recall that the \emph{Ran space}, which we denote by $\Ran$, is the prestack which parameterizes nonempty finite subsets of $X$. It can be written as \[ \Ran = \underset{I}{\colim} \ X^I, \] where the colimit is taken over the category of nonempty finite sets and surjective maps. In particular, a point $S \to \Ran$ factors through a map $x_I : S \to X^I$ for some nonempty finite set $I$.

We denote by $D_{\Ran}$ the multidisk over $\Ran$, defined as the subspace of $X \times \Ran$ consisting of pairs $(x,x_I) : S \to X \times \Ran$ such that \[ \Gamma_x \subset D_{x_I} := \widehat{(X \times S)}_{\Gamma_{x_I}}. \] Here we use the standard notation \[ \Gamma_{x_I} := \bigcup_{i \in I} \Gamma_{x_i} \] for the union of the graphs.

For any prestack $Y$ over $X$ we can form the associated space of arcs $\mathfrak{L}^+_{\Ran}Y$, a prestack over $\Ran$. A point $S \to \mathfrak{L}^+_{\Ran}Y$ consists of $x_I : S \to X^I$ and a map $D_{x_I} \to Y$ over $X$. If $Y$ is constant, i.e. has the form $Z \times X$, we will abuse notation by writing $\mathfrak{L}^+_{\Ran}Z$ rather than $\mathfrak{L}^+_{\Ran}(Z \times X)$.

Let us recall the factorization structure on the arc space. Writing \[ (\Ran \times \Ran)_{\disj} \subset \Ran \times \Ran \] for the set of pairs $(x_I,y_J)$ such that $\Gamma_{x_I} \cap \Gamma_{y_J} = \varnothing$, this structure consists of isomorphisms \[ (\Ran \times \Ran)_{\disj} \times_{\Ran} \mathfrak{L}^+_{\Ran}Y \tilde{\longrightarrow} (\Ran \times \Ran)_{\disj} \times_{\Ran \times \Ran} (\mathfrak{L}^+_{\Ran}Y \times \mathfrak{L}^+_{\Ran}Y), \] where the the map $(\Ran \times \Ran)_{\disj} \to \Ran$ in the first fiber product is given by union of finite subsets of $X$. These isomorphisms satisfy natural associativity and commutativity conditions.

The arc space also admits an evident counital structure compatible with factorization, consisting of a compatible collection of maps \[ \mathfrak{L}^+_{X^I}Y \longrightarrow X^I \times_{X^J} \mathfrak{L}^+_{X^J}Y \] for every injection of (possibly empty) finite sets $J \to I$ (here we put $\mathfrak{L}^+_{X^{\varnothing}}Y := \Spec k$). These maps are isomorphisms in the formal neighborhood of the diagonal (as pointed out in \cite{R} Remark 5.2.2, this is equivalent to the fact that the $\mathfrak{L}^+_{X^I}Y$ assemble into a space over $\Ran$ in the sense of derived algebraic geometry).

\subsection{} Given $x_I : S \to X^I$, we define the associated punctured multidisk by \[ \mathring{D}_{x_I} := \aff(D_{x_I}) \setminus \Gamma_{x_I}. \] In particular, we have maps \[ D_{x_I} \longrightarrow \aff(D_{x_I}) \longleftarrow \mathring{D}_{x_I}. \]

If $Y$ is an affine scheme, we define the associated loop space: a point $S \to \mathfrak{L}_{\Ran}Y$ consists of $x_I : S \to X^I$ along with a map $\mathring{D}_{x_I} \to Y$. Similarly to $\mathfrak{L}^+_{\Ran}Y$, the prestack $\mathfrak{L}_{\Ran}Y$ is equipped with a natural factorization structure. Since \[ \Map(D_{x_I},Y) = \Map(\aff(D_{x_I}),Y), \] we have a canonical embedding of factorization spaces \[ \mathfrak{L}^+_{\Ran}Y \longrightarrow \mathfrak{L}_{\Ran}Y. \]

\subsection{} Observe that the functors $Y \mapsto \mathfrak{L}^+_{\Ran}Y,\mathfrak{L}_{\Ran}Y$ from affine schemes to factorization spaces preserve products. In particular, the factorization spaces $\mathfrak{L}^+_{\Ran}\G_m$ and $\mathfrak{L}_{\Ran}\G_m$ have compatible structures of commutative groups over $\Ran$.

In what follows we will consider the commutative group prestack \[ \Pic(D)_{\Ran} := B_{\Ran}\mathfrak{L}^+\G_m. \] One checks using the smoothness of $\G_m$ that the canonical map \[ \Pic(D)_{\Ran} \longrightarrow \mathfrak{L}^+_{\Ran}B\G_m \] is an isomorphism. That is, a point of $\Pic(D)_{x_I}$ can be viewed as a line bundle on $D_{x_I}$, and any such line bundle is trivial locally on $S$.

We similarly define \[ \Pic(\mathring{D})_{\Ran} := B_{\Ran}\mathfrak{L}_{\Ran}\G_m. \] A point of $\Pic(\mathring{D})_{x_I}$ can be viewed as a line bundle on $\mathring{D}_{x_I}$ which is trivial locally on $S$. It is well-known that the restriction map \[ \Map(\aff(D_{x_I}),B\G_m) \longrightarrow \Map(D_{x_I},B\G_m) \] is an isomorphism (see \cite{Ga1} Lemma 1.1.2), so we can equivalently view $\Pic(\mathring{D})_{x_I}$ as the stack of line bundles on $\mathring{D}_{x_I}$ which, locally on $S$, extend to $D_{x_I}$ (or rather extend to a line bundle on $\aff(D_{x_I})$, which is the same as a line bundle on $D_{x_I}$). In these terms, the homomorphism \[ \Pic(D)_{\Ran} \longrightarrow \Pic(\mathring{D})_{\Ran} \] induced by $\mathfrak{L}^+_{\Ran}\G_m \to \mathfrak{L}_{\Ran}\G_m$ corresponds to restriction of line bundles.

\subsection{} We now introduce the local pairing \[ \langle \, , \, ]_{\mathring{D}_{\Ran}} : \Pic(\mathring{D})_{\Ran} \times_{\Ran} \mathfrak{L}_{\Ran}\G_m \longrightarrow B\G_m. \]

The determinant gerbe construction determines the morphism of group prestacks over $\Ran$
\begin{align*}
\mathfrak{L}_{\Ran}\G_m &\longrightarrow B\G_m \times \Ran \\
(x_I,f) &\mapsto \textstyle (D(f;K_{x_I}),x_I),
\end{align*}
where we view \[ K_{x_I} := \Gamma(\mathring{D}_{x_I},\mathcal{O}) \] as a Tate vector bundle on $S$.

Crucially, this is \emph{not} a morphism of commutative group prestacks. It gives rise to a central extension \[ 1 \longrightarrow \G_m \times \Ran \longrightarrow (\mathfrak{L}_{\Ran}\G_m)^{\flat} \longrightarrow \mathfrak{L}_{\Ran}\G_m \longrightarrow 1, \] which in turn determines the commutator pairing 
\begin{equation}
\label{ccpair}
\mathfrak{L}_{\Ran}\G_m \times_{\Ran} \mathfrak{L}_{\Ran}\G_m \longrightarrow \G_m,
\end{equation}
often referred to as the \emph{Contou-Carr\`{e}re pairing} (see also Section 3 of \cite{BBE} for a similar construction of this pairing). Evidently (\ref{ccpair}) is bilinear and skew-symmetric. In the appendix we show that at a fixed point in $X$, this pairing agrees with the one actually introduced by Contou-Carr\`{e}re in \cite{CC1}.

Finally, we define the local pairing \[ \langle \, , \, ]_{\mathring{D}_{\Ran}} : \Pic(\mathring{D})_{\Ran} \times_{\Ran} \mathfrak{L}_{\Ran}\G_m \longrightarrow B\G_m \] by delooping (\ref{ccpair}) in the first variable.

\begin{proposition}
The pairing  $\langle \, , \, ]_{\mathring{D}_{\Ran}}$ is given by the formula \[ \textstyle \langle \mathcal{L} , f ]_{\mathring{D}_{\Ran}} = D(f;\Gamma(\mathring{D}_{x_I},\mathcal{L}))^{-1} \otimes D(f;K_{x_I}) \] where $\mathcal{L}$ is an extendable line bundle on $\mathring{D}_{x_I}$ and $f : \mathring{D}_{x_I} \to \G_m$.
\end{proposition}

\begin{proof}
The displayed formula determines a homomorphism 
\begin{equation}
\label{ccform}
\mathfrak{L}_{\Ran}\G_m \longrightarrow \uMap_{\Ran,\bullet}(\Pic(\mathring{D})_{\Ran},B\G_m \times \Ran) \cong \uHom_{\Ran}(\mathfrak{L}_{\Ran}\G_m,\G_m \times \Ran),
\end{equation}
which we must show agrees with (\ref{ccpair}).

Applying the pairing (\ref{ccpair}) to $f,g : \mathring{D}_{x_I} \to \G_m$ yields the automorphism
\begin{align*}
\textstyle D(fg,K_{x_I}) &\tilde{\longrightarrow} \textstyle D(f,K_{x_I}) \otimes D(g,K_{x_I}) \\
&\tilde{\longrightarrow} \textstyle D(g,K_{x_I}) \otimes D(f,K_{x_I}) \\
&\tilde{\longrightarrow} \textstyle D(gf,K_{x_I}) = D(fg,K_{x_I}),
\end{align*}
interpreted as an invertible function on $S$.

On the other hand, the pairing (\ref{ccform}) applied to $f,g$ gives the invertible function on $S$ corresponding to the automorphism \[ D(f,K_{x_I}) \overset{g}{\tilde{\longrightarrow}} D(gfg^{-1},K_{x_I}) = D(f,K_{x_I}). \] By inspection, these functions agree.
\end{proof}

\subsection{} Let us recall the definition of the Beilinson-Drinfeld affine Grassmannian $\Gr(D)_{\Ran}$ for $\G_m$. An $S$-point of $\Gr(D)_{\Ran}$ consists of $x_I : S \to X^I$, a line bundle $\mathcal{L}$ on $X \times S$, and a trivialization of $\mathcal{L}$ over $(X \times S) \setminus \Gamma_{x_I}$.

By the Beauville-Laszlo theorem, a point of $\Gr(D)_{x_I}$ can equivalently be described as a line bundle on $D_{x_I}$ together with a trivialization over $\mathring{D}_{x_I}$. In particular, we have a short exact sequence \[ 1 \longrightarrow \mathfrak{L}^+_{\Ran}\G_m \longrightarrow \mathfrak{L}_{\Ran}\G_m \longrightarrow \Gr(D)_{\Ran} \longrightarrow 1 \] of fppf sheaves relative to $\Ran$ (i.e. the fibers over any $x_I : S \to X^I$ form a short exact sequence of fppf sheaves of abelian groups over $S$).

As is well-known, the space $\Gr(D)_{\Ran}$ admits a natural factorization structure. Moreover $\Gr(D)_{\Ran}$ admits a unital structure compatible with factorization. This structure consists of a compatible system of maps \[ X^I \times_{X^J} \Gr(D)_{X^J} \longrightarrow \Gr(D)_{X^I} \] over $X^I$ for every injection of (possibly empty) finite sets $J \to I$.

\subsection{} We also consider the submonoid $\Gr_+(D)_{\Ran} \subset \Gr(D)_{\Ran}$ consisting of points where the given section extends to $X \times S$, possibly with zeros.

\begin{proposition}
\label{grrangrp}
The inclusion $\Gr_+(D)_{\Ran} \to \Gr(D)_{\Ran}$ realizes the latter as the group completion of the former over $\Ran$.
\end{proposition}

\begin{proof}
Fix $x_I : S \to X^I$, a line bundle $\mathcal{L}$ on $X \times S$, and a nonvanishing section $\sigma$ of $\mathcal{L}$ defined away from $\Gamma_{x_I}$. Viewing $\Gamma_{x_I}$ as an effective Cartier divisor, we see that \[ \sigma \in \Gamma(X \times S,\mathcal{L}(n \cdot \Gamma_{x_I})) \] for $n \gg 0$. Thus we obtain the required factorization \[ (\mathcal{L},\sigma) = (\mathcal{L}(n \cdot \Gamma_{x_I}),\sigma) \otimes (\mathcal{O}_{X \times S}(n \cdot \Gamma_{x_I}),1)^{-1}. \]
\end{proof}

Denote by $\Div_+(D)_{\Ran}$ the space of effective relative Cartier divisors on $D_{\Ran}$, a sheaf of commutative monoids relative to $\Ran$. A point $S \to \Div_+(D)_{\Ran}$ consists of a point $x_I : S \to X^I$ of $\Ran$ together with an effective divisor on $X \times S$ relative to $S$ which is contained in $D_{x_I}$. Note that $\Div_+(D)_{X^I}$ is a scheme for all $I$, with connected components $\Div^d_+(D)_{X^I}$ labeled by $d \geq 0$. It is generally highly reducible and nonreduced.

We have a tautological isomorphism
\begin{align}
\label{divgrloc}
\Div_+(D)_{\Ran} &\tilde{\longrightarrow} \Gr_+(D)_{\Ran} \\
(x_I,Z) &\mapsto (x_I,\mathcal{O}_{X \times S}(Z),1). \nonumber
\end{align}
In what follows, we will use this identification without comment and view points of $\Gr(D)_{\Ran}$ as Cartier divisors when convenient. We write $\Gr^d_+(D)_{X^I} := \Div^d_+(D)_{X^I}$, etc.

\subsection{} Given a line bundle $\mathcal{L}$ on $D_{x_I}$, note that $\Gamma(D_{x_I},\mathcal{L})$ is a lattice in the Tate vector bundle $\Gamma(\mathring{D}_{x_I},\mathcal{L})$ (the claim is local on $S$, hence reduces to the case $\mathcal{L} = \mathcal{O}$).

\begin{proposition}
\label{unrpairing}
The restricted pairing \[ \Pic(D)_{\Ran} \times_{\Ran} \mathfrak{L}^+_{\Ran}\G_m \longrightarrow \Pic(\mathring{D})_{\Ran} \times_{\Ran} \mathfrak{L}_{\Ran}\G_m \xrightarrow{\langle \, , \, ]_{\mathring{D}_{\Ran}}} B\G_m \] is canonically trivial, hence determines a pairing \[ \langle \, , \, ]_{D_{\Ran}} : \Pic(D)_{\Ran} \times_{\Ran} \Gr(D)_{\Ran} \longrightarrow B\G_m. \] This pairing is given by the formula \[ \langle \mathcal{L},Z]_{D_{\Ran}} = \det(\Gamma(D_{x_I},\mathcal{L}),\Gamma(D_{x_I},\mathcal{L}(Z)))^{-1} \otimes \det(\Gamma(D_{x_I},\mathcal{O}),\Gamma(D_{x_I},\mathcal{O}(Z))), \] or when $Z$ is effective by \[ \langle \mathcal{L},Z]_{D_{\Ran}} = \det \Gamma(Z,\mathcal{L}) \otimes (\det \Gamma(Z,\mathcal{O}))^{-1}. \]
\end{proposition}

\begin{proof}
Fix $x_I : S \to X^I$, a line bundle $\mathcal{L}$ on $D_{x_I}$, and $f : D_{x_I} \to \G_m$. Then the line \[ D(f;\Gamma(\mathring{D}_{x_I},\mathcal{L})) = \det(f\Gamma(D_{x_I},\mathcal{L}),\Gamma(D_{x_I},\mathcal{L})) = \det(\Gamma(D_{x_I},\mathcal{L}),\Gamma(D_{x_I},\mathcal{L})) \] is canonically trivial. The previous proposition yields a trivialization of the line $\langle \mathcal{L}|_{\mathring{D}_{x_I}},f ]_{\mathring{D}_{\Ran}}$.

The formulas in question follow from the observation that for $Z = \ddiv(f)$, we have \[ f\Gamma(D_{x_I},\mathcal{L}) = \Gamma(D_{x_I},\mathcal{L}(-Z)) \] as lattices in $\Gamma(\mathring{D}_{x_I},\mathcal{L})$.
\end{proof}

\section{The local Albanese property}

\label{localbsec}

\subsection{} If $G$ and $H$ are two fppf sheaves of abelian groups over $S$, we write $\uExt_S^1(G,H)$ for the sheaf over $S$ associated to the presheaf
\[ T \mapsto \Ext^1_T(G \times_S T, H \times_S T). \]

The next lemma, which seems to be new, is an important input into Theorem \ref{locramthm}.

\begin{lemma}
\label{extvanish}
Let $G$ be a commutative affine group scheme over $S$ such that $\mathcal{O}_G$ is a projective $\mathcal{O}_S$-module. Then the sheaf $\uExt^1_S(G, \G_m)$ vanishes.
\end{lemma}

\begin{proof}
We need to show that any extension of fppf sheaves
\[ 0 \longrightarrow \G_m \times T \longrightarrow E \longrightarrow G \times_S T \longrightarrow 0 \]
splits locally on $T$. Since $E \cong \G_m \times G \times_S T$ locally on $G \times_S T$, we see that $\mathcal{O}_E$ is a projective module over $\mathcal{O}_G \otimes_{\mathcal{O}_S} \mathcal{O}_T$ (recall that projectivity is fppf local). Since $\mathcal{O}_G$ is projective over $\mathcal{O}_S$, it follows that $\mathcal{O}_E$ is a projective $\mathcal{O}_T$-module.

We now apply Theorem 3.5 in \cite{A} to see that $E$ and $G$ are $0$-reflexive as sheaves of abelian groups over $T$. Hence the $0$-dual sequence 
\[ 0 \longrightarrow G^* \times_S T \longrightarrow E^* \longrightarrow \Z \times T \longrightarrow 0 \]
is exact. The fiber of $E^* \to \Z \times T$ over $\{1\} \times T$ is therefore a locally trivial torsor on $T$, and a local section of this torsor determines a local splitting of $E^*$. Dualizing again, we obtain a local splitting of $E$.
\end{proof}

\begin{corollary} 
Any commutative affine group scheme over $S$ as in the lemma is $1$-reflexive.
\end{corollary}

\begin{proof}
This follows from the previous lemma combined with Corollary $3.6$ in \cite{Br}.
\end{proof}

In the factorization setting, we will apply the previous lemma to the arc group.

\begin{lemma}
\label{arcprojlem}
For any nonempty finite set $I$, the structure sheaf of $\mathfrak{L}^+_{X^I}\G_m$ is locally projective as an $\mathcal{O}_{X^I}$-module.
\end{lemma}
\begin{proof}
Observe that as a factorization space, the arc group is universal in the sense of \cite{Cl}. It therefore suffices to treat the case $X = \mathbb{A}^1$. We proceed by induction on the cardinality of $I$. 

For the base case, note that by choosing a global coordinate on $X$, we obtain an isomorphism
\[ \mathfrak{L}^+_X \G_m \longrightarrow \mathfrak{L}^+_0 \G_m \times X \]
over $X$, where $\mathfrak{L}^+_0 \G_m$ denotes the fiber of $\mathfrak{L}^+_X \G_m$ over $0$. Now assume the claim for all finite sets $J$ with $|J| < |I|$ and let $\Delta \subseteq X^I$ denote the main diagonal. Since local projectivity of a sheaf is an fpqc local property, it suffices to check the claim on the cover of $X^I$ consisting of the complement of the diagonal and the affinization of the formal completion of $X^I$ along $\Delta$. Note that this cover is indeed fpqc since $\mathcal{O}_{X^I}$ is Noetherian, and so the natural map from $\mathcal{O}_{X^I}$ to its completion is flat. 

Away from $\Delta$, we arrive at our desired conclusion using the factorization structure on the arc group and our induction hypothesis. It therefore suffices to check that the structure sheaf of $\mathfrak{L}^+_{X^I}\G_m$ is locally projective when restricted to the affinization of the formal neighborhood of $\Delta$. Recall that the counital structure on the arc group provides us, for every injection $J \to I$, with morphisms 
\[\mathfrak{L}^+_{X^I} \G_m \longrightarrow X^I \times_{X^J} \mathfrak{L}^+_{X^J} \G_m, \]
which are isomorphisms along the formal neighborhood of $\Delta$. In particular, we may take $J$ to be a singleton, reducing to our base case. Hence we see that the structure sheaf of $\mathfrak{L}^+_{X^I}\G_m$ is in fact free when restricted to the affinized completion of $X^I$ along $\Delta$. 
\end{proof}

\subsection{} Consider the Abel-Jacobi map over $\Ran$
\begin{align*}
\AJ_{D_{\Ran}} : D_{\Ran} &\longrightarrow \Gr(D)_{\Ran} \\
    (x_I,y) &\mapsto (x_I,\mathcal{O}_{X \times S}(y),1) \nonumber
\end{align*}
attached to the multidisk $D_{\Ran}$. It is an isomorphism onto $\Gr^1_+(D)_{\Ran}$.

Notice that the unital factorization structure on $\Gr(D)_{\Ran}$ induces a counital factorization structure on its Cartier dual $\Gr(D)^{\vee}_{\Ran}$ taken over $\Ran$, compatibly with its commutative group structure.

We now prove the unramified local Cartier duality theorem. Recall the pairing $\langle \, , \, ]_{D_{\Ran}}$ from Proposition \ref{unrpairing}.

\begin{theorem}
\label{locunrthm}
The pairing $\langle \, , \, ]_{D_{\Ran}}$ is perfect. The resulting isomorphism \[ \Gr(D)^{\vee}_{\Ran} \tilde{\longrightarrow} \Pic(D)_{\Ran} \] is an isomorphism of counital factorization spaces, and agrees with restriction along $\AJ_{D_{\Ran}}$.
\end{theorem}

\begin{proof}[Proof]

Let us verify that restriction along $\AJ_{D_{\Ran}}$ is inverse to the homomorphism  \[ \Pic(D)_{\Ran} \longrightarrow \Gr(D)_{\Ran}^{\vee} \] induced by $\langle \, , \, ]_{D_{\Ran}}$. Given a section $y : S \to D_{x_I}$, we have a tautological isomorphism $\langle \mathcal{L},\Gamma_y]_{D_{\Ran}} \cong y^*\mathcal{L}$, which shows that \[ \Pic(D)_{\Ran} \longrightarrow \Gr(D)^{\vee}_{\Ran} \longrightarrow \Pic(D)_{\Ran} \] is isomorphic to the identity.

Conversely, by Proposition \ref{grrangrp} it suffices to show that for any multiplicative line bundle $\mathcal{M}$ on $\Gr_{x_I,+}$, a trivialization of $\mathcal{M}|_{D_{x_I}} = \AJ_{D_{x_I}}^*\mathcal{M}$ determines a trivialization of $\mathcal{M}$. By induction, we may assume that we have constructed a trivialization of $\mathcal{M}$ over $\Gr_{x_I,+}^d$ for all $d \leq n$.

Observe that the addition morphism \[ \add : D_{x_I} \times_S \Gr_{x_I,+}^n \longrightarrow \Gr_{x_I,+}^{n+1} \] canonically identifies with the universal effective divisor of degree $n+1$, and in particular this morphism is a finite flat covering. The inductive hypothesis yields a trivialization of \[ \add^*\mathcal{M} \cong \mathcal{M}|_{D_{x_I}} \boxtimes_{\mathcal{O}_S} \mathcal{M}|_{\Gr_{x_I,+}^n}, \] and the cocycle condition on the multiplicative structure of $\mathcal{M}$ ensures that this trivialization descends to $\Gr_{x_I,+}^{n+1}$.

We see by inspection that restriction along $\AJ_{D_{\Ran}}$  is compatible with counital factorization structures.

It remains to show that $\Pic(D)_{\Ran}$ is $1$-reflexive, from which it will follow that $\Gr(D)_{\Ran}$ is $1$-reflexive. Our argument thus far implies that restriction along $\AJ_{D_{\Ran}}$ induces an isomorphism \[ \Gr(D)_{\Ran}^* = H^{-1}(\Gr(D)_{\Ran}^{\vee}) \tilde{\longrightarrow} H^{-1}(\Pic(D)_{\Ran}) = \mathfrak{L}^+_{\Ran}\G_m, \] and that $\uExt^1_{\Ran}(\Gr(D)_{\Ran},\G_m) = 1$. By Lemma \ref{arcprojlem} and Theorem 3.5 in \cite{A}, the arc group $\mathfrak{L}^+_{\Ran}\G_m$ is $0$-reflexive, and hence \[ (\mathfrak{L}^+_{\Ran}\G_m)^* \tilde{\longrightarrow} \Gr(D)_{\Ran}. \] Since $\Pic(D)_{\Ran} = B_{\Ran}\mathfrak{L}^+_{\Ran}\G_m$, we have \[ \Pic(D)_{\Ran}^{\vee} = (\mathfrak{L}^+_{\Ran}\G_m)^*. \] Dualizing again, we obtain \[ \Pic(D)_{\Ran}^{\vee \vee} = B(\mathfrak{L}^+_{\Ran}\G_m)^{**}, \] so we are done by the $0$-reflexivity of $\mathfrak{L}^+_{\Ran}\G_m$.
\end{proof}

\begin{corollary}
\label{locunralb}
For any $1$-reflexive group stack $G$ over $\Ran$, restriction along $\AJ_{D_{\Ran}}$ induces an isomorphism \[ \uHom_{\Ran}(\Gr(D)_{\Ran},G) \longrightarrow \uMap_{\Ran}(D_{\Ran},G). \]
\end{corollary}

\begin{proof}
We have
\begin{align*}
\uHom_{\Ran}(\Gr(D)_{\Ran},G) &\tilde{\longrightarrow} \uHom_{\Ran}(G^{\vee},\Gr(D)^{\vee}_{\Ran}) \\
&\tilde{\longrightarrow} \uHom_{\Ran}(G^{\vee},\Pic(D)_{\Ran}) \\
&\tilde{\longrightarrow} \uMap_{\Ran}(D_{\Ran},G^{\vee \vee}) \\
&\tilde{\longrightarrow} \uMap_{\Ran}(D_{\Ran},G),
\end{align*}
which is easily seen to agree with restriction along $\AJ_{D_{\Ran}}$.
\end{proof}

\subsection{} The next result says that the forgetful map \[ \Gr(D)_{\Ran} \longrightarrow \Pic(D)_{\Ran} \] is self-dual under the pairing $\langle \, , \, ]_{D_{\Ran}}$.

\begin{proposition}
\label{locpairingcomm}
The square \[
\begin{tikzcd}
\Gr(D)_{\Ran} \times_{\Ran} \Gr(D)_{\Ran} \arrow{r} \arrow{d} & \Pic(D)_{\Ran} \times_{\Ran} \Gr(D)_{\Ran} \arrow{d}{\langle \, , \, ]_{D_{\Ran}}} \\
\Gr(D)_{\Ran} \times_{\Ran} \Pic(D)_{\Ran} \arrow{r}{[ \, , \, \rangle_{D_{\Ran}}} & B\G_m
\end{tikzcd} \]
commutes up to canonical isomorphism.
\end{proposition}

\begin{proof}

By Corollary \ref{locunralb}, it suffices to construct an isomorphism between the two pairings after restriction along \[ \AJ_{D_{\Ran}} \times_{\Ran} \AJ_{D_{\Ran}} : D_{\Ran} \times_{\Ran} D_{\Ran} \longrightarrow \Gr(D)_{\Ran} \times_{\Ran} \Gr(D)_{\Ran}. \] Observe that for any $x_I : S \to X^I$ and any $y,z : S \to D_{x_I}$, we have a canonical isomorphism \[ \langle \mathcal{O}_{D_{x_I}}(y),\Gamma_z]_{D_{\Ran}} \tilde{\longrightarrow} (z,y)^*\mathcal{O}_{D_{x_I} \times_S D_{x_I}}(\Delta), \] where $\Delta : D_{x_I} \to D_{x_I} \times_S D_{x_I}$ denotes the diagonal. Now the claim follows from the equivariance of $\Delta$ under transposition of the factors.

\end{proof}

\begin{lemma}
\label{loopextvanish}
We have $\uExt^1_{\Ran}(\mathfrak{L}_{\Ran}\G_m,\G_m) = 1$.
\end{lemma}
\begin{proof}
Combining Lemmas \ref{extvanish} and \ref{arcprojlem}, we obtain $\uExt^1_{\Ran}(\mathfrak{L}^+_{\Ran}\G_m,\G_m) = 1$. On the other hand $\uExt^1_{\Ran}(\Gr(D)_{\Ran},\G_m) = 1$ by Theorem \ref{locunrthm}, whence the lemma.
\end{proof}

\subsection{} We are now in a position to prove the local Cartier duality theorem.

\begin{theorem}
\label{locnondeg}
The pairing $\langle \, , \, ]_{\mathring{D}_{\Ran}}$ is perfect.
\end{theorem}

\begin{proof}
Consider the diagram \[
\begin{tikzcd}
\Gr(D)_{\Ran} \arrow{r} \arrow{d} & \Pic(D)_{\Ran} \arrow{r} \arrow{d} & \Pic(\mathring{D})_{\Ran} \arrow{d} \\
\Pic(D)_{\Ran}^{\vee} \arrow{r} & \Gr(D)_{\Ran}^{\vee} \arrow{r} & \mathfrak{L}_{\Ran}\G_m^{\vee}.
\end{tikzcd} \]
Here the duals are taken over $\Ran$, the left and center vertical maps are induced by $\langle \, , \, ]_{D_{\Ran}}$, and the right vertical map is induced by $\langle \, , \, ]_{\mathring{D}_{\Ran}}$. Note that the bottom row is exact by Lemma \ref{loopextvanish}.

The right square commutes by the construction of $\langle \, , \, ]_{D_{\Ran}}$, and the left square commutes by Proposition \ref{locpairingcomm}. Since the left and center vertical maps are isomorphisms by Theorem \ref{locunrthm}, the right vertical map is an isomorphism as desired.
\end{proof}

Having established the purely local part of Theorem \ref{locramthm}, we can already deduce the local Albanese property, i.e. Corollary \ref{locramalb}. Let $G$ be a $1$-reflexive commutative group stack relative to $\Ran$, meaning $G$ is a commutative group prestack over $\Ran$ such that for any $S \to \Ran$, the fiber product $G \times_{\Ran} S$ is a $1$-reflexive commutative group stack over $S$. Then we define the stack of $B\G_m$-extendable loops \[ \mathfrak{L}^{\ext}_{\Ran}G := \uHom_{\Ran}(G^{\vee},\Pic(\mathring{D})_{\Ran}). \] We have a fully faithful embedding \[ \mathfrak{L}^{\ext}_{\Ran}G \longrightarrow \mathfrak{L}_{\Ran}G \] whose image can be characterized analogously to what was done in Section \ref{extdefsec}, replacing $X$ by $\mathring{D}$ and $\overline{X}$ by $D$.

We now formulate and prove a more precise version of Corollary \ref{locramalb}.

\begin{corollary}
\label{locramalb2}
There is a unique isomorphism \[ \uHom_{\Ran}(\mathfrak{L}_{\Ran}\G_m,G) \tilde{\longrightarrow} \mathfrak{L}^{\ext}_{\Ran}G \] which is natural in $1$-reflexive commutative group stacks $G$ relative to $\Ran$ and specializes to the isomorphism of Theorem \ref{locnondeg} when $G = B\G_m$. Moreover, this isomorphism is compatible with the factorization structures on both sides.
\end{corollary}

\begin{proof}
Compose the isomorphisms
\begin{align*}
\uHom_{\Ran}(\mathfrak{L}_{\Ran}\G_m,G) &\tilde{\longrightarrow} \uHom_{\Ran}(G^{\vee},\mathfrak{L}_{\Ran}\G_m^{\vee}) \\
&\tilde{\longrightarrow} \uHom_{\Ran}(G^{\vee},\Pic(\mathring{D})_{\Ran}) = \mathfrak{L}^{\ext}_{\Ran}G.
\end{align*}
The compatibility with factorization structures is manifest.
\end{proof}

\subsection{} In some cases of interest, the $B\G_m$-extendability condition is unnecessary. For simplicity we work over a fixed $x \in X(k)$ here.

\begin{proposition}
Suppose that $G$ is a commutative affine group scheme over $k$. Then we have a canonical isomorphism \[ \mathfrak{L}_x^{\ext}G \longrightarrow \mathfrak{L}_xG. \]
\end{proposition}

\begin{proof}
By Lemma \ref{extvanish} we have $G^{\vee} \cong BG^*$. Thus \[ \mathfrak{L}_x^{\ext}G = \uHom(G^{\vee},\Pic(\mathring{D}_x)) \cong \uHom(G^*,\mathfrak{L}_x\G_m). \] Since $G$ is affine its $0$-Cartier dual $G^*$ is ind-finite, which implies that \[ \uHom(G^*,\mathfrak{L}_x\G_m) \cong \mathfrak{L}_xG \] as desired.
\end{proof}

Combining the previous proposition with Corollary \ref{locramalb}, we obtain the following result.

\begin{corollary}
If $G$ is as in the proposition, then the isomorphism of Corolllary \ref{locramalb2} specializes to \[ \uHom(\mathfrak{L}_x\G_m,G) \tilde{\longrightarrow} \mathfrak{L}_xG. \]
\end{corollary}

\section{The unramified global Albanese property}

\subsection{} For the entirety of this section, we assume that $X = \overline{X}$ is projective.

\begin{definition}
An open subscheme $U \subset X \times S$ is called a \emph{domain} relative to $S$ if the projection $U \to S$ is surjective.
\end{definition}

We define the global affine Grassmannian $\Gr(X)$ to be the following commutative group stack: a point $S \to \Gr(X)$ consists of a line bundle $\mathcal{L}$ on $X \times S$, a domain $U \subset X \times S$, and a nonvanishing section of $\mathcal{L}$ over $U$. An isomorphism between two $S$-points of $\Gr(X)$ is an isomorphism of line bundles under which the two sections agree on a common domain. This is easily seen to be an equivalence relation, i.e. $\Gr(X)$ is isomorphic to a sheaf of abelian groups.

We will also consider the submonoid $\Gr_+(X) \subset \Gr(X)$ consisting of those points such that the given section extends to $X \times S$ (possibly with zeros).

Note that $(X \times S) \setminus \Gamma_{x_I}$ is a domain relative to $S$ for any $x_I : S \to X^I$. This implies the existence of a canonical map
\begin{equation}
\label{rantorat}
\Gr(D)_{\Ran} \longrightarrow \Gr(X)
\end{equation}
which sends $\Gr_+(D)_{\Ran}$ into $\Gr_+(X)$.

\begin{proposition}
The map (\ref{rantorat}) is an fppf local surjection.
\end{proposition}

\begin{proof}
A standard argument shows that for any domain $U \subset X \times S$, fppf locally on $S$ there exists $x_I : S \to X^I$ such that $(X \times S) \setminus \Gamma_{x_I} \subset U$ (see \cite{Lu}, Proposition 2). 
\end{proof}

\begin{corollary}
\label{grratgrp}
The inclusion $\Gr_+(X) \to \Gr(X)$ realizes the latter as the group completion of the former in the category of fppf sheaves.
\end{corollary}
\begin{proof}
Combine the previous proposition with Proposition \ref{grrangrp}.
\end{proof}

Let $\Div_+(X)$ denote the moduli of effective Cartier divisors on $X$, which is a sheaf of commutative monoids. It is well-known that the map
\begin{align*}
\coprod_{d \geq 0} X^d &\longrightarrow \Div_+(X) \\
(x_i) &\mapsto \sum_i x_i
\end{align*}
factors through an isomorphism \[ \coprod_{d \geq 0} X^{(d)} \tilde{\longrightarrow} \Div_+(X), \] where $X^{(d)} := X^d/\!/\Sigma_d$ is the $d^{\text{th}}$ symmetric power. Moreover, we have the evident isomorphism
\begin{align*}
\Div_+(X) &\tilde{\longrightarrow} \Gr_+(X) \\
Z &\mapsto (\mathcal{O}_X(Z),1).
\end{align*}

\subsection{} Observe that the map $\AJ_X$ factors through
\begin{align}
\label{ajgr}
X &\longrightarrow \Gr(X) \\
x &\mapsto (\mathcal{O}_X(x),1). \nonumber
\end{align}

\begin{lemma}
\label{globgrdual}
Restriction along (\ref{ajgr}) induces an isomorphism \[ \Gr(X)^{\vee} \tilde{\longrightarrow} \Pic(X). \]
\end{lemma}

The lemma will be deduced from Theorem \ref{locunrthm} by local-to-global methods. Recall that for any prestack $\mathcal{Z}$, we can consider its de Rham prestack $\mathcal{Z}_{\dR}$, defined by \[ \Map(S,\mathcal{Z}_{\dR}) = \Map(S_{\red},\mathcal{Z}). \]

We will also need to consider the space $\mathfrak{L}_{\Ran}^{+,\hor}Y$ of horizontal jets into a prestack $Y \to X_{\dR}$. See Section 5 of \cite{R}, where the construction is called ``multijets," for more details. A point $S \to \mathfrak{L}_{\Ran}^{+,\hor}Y$ consists of $x_I : S \to X^I$ together with a map \[ (\Gamma_{x_I})_{\dR} \times_{S_{\dR}} S \longrightarrow Y \] over $X_{\dR}$. It is not difficult to see that $\mathfrak{L}_{\Ran}^{+,\hor}Y$ has a natural structure of counital factorization space, and moreover we have \[ \mathfrak{L}_X^{+,\hor}Y = Y \times_{X_{\dR}} X. \] Restriction along the map \[ D_{x_I} \longrightarrow (D_{x_I})_{\dR} \times_{S_{\dR}} S = (\Gamma_{x_I})_{\dR} \times_{S_{\dR}} S \] defines a morphism of counital factorization spaces \[ \mathfrak{L}_{\Ran}^{+,\hor}Y \longrightarrow \mathfrak{L}_{\Ran}^+(Y \times_{X_{\dR}} X). \]

\begin{proof}[Proof of Lemma \ref{globgrdual}]
Let $\Pic(D)_{X_{\dR}}$ denote the space over $X_{\dR}$ defined as follows: a point $S \to \Pic(D)_{X_{\dR}}$ consists of $x : S_{\red} \to X$ together with a line bundle on $D_x$ (note that the completion along $\Gamma_x$ is well-defined even though $x$ is only defined on $S_{\red}$). Clearly we have \[ \Pic(D)_X = \Pic(D)_{X_{dR}} \times_{X_{\dR}} X. \]

Tracing through the definitions, we have a tautological isomorphism of counital factorization spaces \[ \Pic(D)_{\Ran} \tilde{\longrightarrow} \mathfrak{L}_{\Ran}^{+,\hor}(\Pic(D)_{X_{\dR}}). \] Consider the moduli stack of sections of the projection $\Pic(D)_{\Ran} \to \Ran$ compatible with the counital factorization structure. By Theorem 5.2.3 in \cite{R}, this stack of factorizable sections identifies with the stack of sections of $\Pic(D)_{X_{\dR}} \to X_{\dR}$. The latter is tautologically isomorphic to $\Pic(X)$.

We claim that the map $\Gr(D)_{\Ran} \to \Gr(X)$ induces an isomorphism from $\Gr(X)^{\vee}$ to the moduli stack of sections of $\Gr(D)^{\vee}_{\Ran} \to \Ran$ compatible with the counital factorization structure. Then the lemma will follow from Theorem \ref{locunrthm} by commutativity of the square \[
\begin{tikzcd}
D_{\Ran} \arrow{r}{\AJ_{D_{\Ran}}} \arrow{d} & \Gr(D)_{\Ran} \arrow{d} \\
X \arrow{r} & \Gr(X).
\end{tikzcd} \]

It follows from Proposition 5.2.2 of \cite{Bar} that the canonical map \[ |\Ran^{\bullet} \times \Gr(D)_{\Ran}| \longrightarrow \Gr(X) \] is an isomorphism locally in the fppf topology, where the left side denotes the geometric realization of the simplicial prestack determined by the action of $\Ran$ on $\Gr(D)_{\Ran}$ coming from the unital structure. One checks from the definitions that multiplicative line bundles on this geometric realization are precisely counital factorizable sections of $\Gr(D)^{\vee}_{\Ran}$, as desired.
\end{proof}

\subsection{} A standard argument shows that any family of line bundles on $X$ trivializes over a domain, i.e. the projection $\Gr(X) \to \Pic(X)$ is an fppf local surjection. Thus we have an exact triangle of commutative group stacks \[ \Rat(X) \longrightarrow \Gr(X) \longrightarrow \Pic(X). \] Here $\Rat$ is the moduli of invertible rational functions on $X$.

We will need the following lemma.

\begin{lemma}
\label{extend}
Let $V$ be a finite-dimensional vector space and $U \subset \mathbb{P}(V)$ an open subscheme whose complement has codimension at least $3$. Then the restriction map \[ \uMap(\mathbb{P}(V),B\G_m) \longrightarrow \uMap(U,B\G_m) \] is an isomorphism.
\end{lemma}

\begin{proof}
Fix a test scheme $S$. Extending scalars if necessary, we may assume that $U$ has a rational point $p \in U(k)$. Denote by $C$ the projective cone on $\mathbb{P}(V) \setminus U$ with vertex $p$, and put $U_0 := \mathbb{P}(V) \setminus C$. Since $C$ has codimension at least $2$, functions on $U_0 \times S$ extend to $\mathbb{P}(V) \times S$, so the restriction functors \[ \Map(\mathbb{P}(V) \times S,B\G_m) \longrightarrow \Map(U \times S,B\G_m) \longrightarrow \Map(U_0 \times S,B\G_m) \] are fully faithful. Therefore it suffices to show that for any line bundle $\mathcal{L}$ on $U \times S$, the restriction $\mathcal{L}|_{U_0 \times S}$ extends to $\mathbb{P}(V) \times S$.

Let $E \subset U \times U_0$ be the subscheme of pairs $(x,y)$ such that $x$, $y$, and $p$ are collinear. Write $\pi : E \to U$ and $\rho : E \to U_0$ for the projections. Since $\rho$ is a $\mathbb{P}^1$-bundle, the degree of $(\pi \times \id_S)^*\mathcal{L}$ is a locally constant function on $U_0 \times S$. Thus, by twisting $\mathcal{L}$ appropriately we can assume that $(\pi \times \id_S)^*\mathcal{L} \cong (\rho \times \id_S)^*\mathcal{L}'$ for some line bundle $\mathcal{L}'$ on $U_0 \times S$, and we are trying to prove that $\mathcal{L}|_{U_0 \times S}$ is pulled back from $S$.

Restricting the chosen isomorphism along $\Delta_{U_0} \times \id_S$, we obtain $\mathcal{L}|_{U_0 \times S} \cong \mathcal{L}'$. On the other hand, we can restrict the same isomorphism along the section $p \times \id_{U_0} \times \id_S$ of $\rho \times \id_S$, which yields an identification $\mathcal{L}|_{U_0 \times S} \cong \mathcal{O}_{U_0} \boxtimes (\mathcal{L}|_{p \times \id_S})$ as needed.
\end{proof}

\subsection{} Let $K$ denote the $k$-vector space of rational functions on $X$. We present $K$ as an ind-scheme using the formula \[ \underset{Z}{\colim} \ \Gamma(X,\mathcal{O}_X(Z)) = K \] where $Z$ runs through all effective divisors on $X$.

Observe that the ind-scheme $K^{\times} := K \setminus \{ 0 \}$ is a commutative monoid under multiplication, although not a group. There is an injective monoid homomorphism $K^{\times} \to \Rat$, which can be seen as follows. A morphism $S \to K^{\times}$ lands in some $\Gamma(X,\mathcal{O}(Z)) \setminus \{ 0 \}$, hence determines a function $f : (X \setminus Z) \times S \to \mathbb{A}^1$ which does not vanish identically on any fiber over a point of $S$. This means that the complement of the vanishing locus of $f$ is a domain.

\begin{lemma}

The inclusion $K^{\times} \to \Rat(X)$ realizes the latter as the group completion as the former in the category of Zariski sheaves.

\end{lemma}

\begin{proof}

The group completion of $K^{\times}$ injects into $\Rat$ because the multiplication in $K^{\times}$ is cancellative. Thus it remains to prove surjectivity.

Fix a domain $U \subset X \times S$ and $f : U \to \G_m$. Choose a closed point $x$ in $X$ and put $\mathring{X} := X \setminus \{ x \}$. Shrink $U$ if necessary so that $U \subset \mathring{X} \times S$, and fix a closed point $s$ in $S$. Then there is a function $g$ on $\mathring{X} \times S$ which vanishes on the complement of $U$ but does not vanish identically on the fiber $\mathring{X}_s$. Replacing $S$ by a Zariski neighborhood of $s$, we may assume $g$ does not vanish identically on any fiber over a point of $S$. Now there exists $m$ so that $fg^m$ is regular on $\mathring{X} \times S$, which means for sufficiently large $n$ we have \[ fg^m \in \Gamma(X \times S,\mathcal{O}(nx) \boxtimes \mathcal{O}_S) \] and likewise for $g^m$. Thus $f$ can be written as the ratio of two $S$-points of $\Gamma(X,\mathcal{O}(nx)) \setminus \{ 0 \} \subset K^{\times}$.

\end{proof}

The previous lemma implies that $K^{\times} \times K^{\times} \to \Rat(X)$ given by $(f,g) \mapsto f/g$ is a Zariski epimorphism. Moreover, this map is evidently $\G_m$-equivariant with respect to the diagonal action on $K^{\times} \times K^{\times}$ and the trivial action on $\Rat(X)$, so we obtain
\begin{equation}
(K^{\times} \times K^{\times})/\G_m \longrightarrow \Rat(X).
\label{ratgrp}
\end{equation}
Observe that the domain is the ind-projective space $\mathbb{P}(K \oplus K)$ with the subspaces $\mathbb{P}(K \oplus 0)$ and $\mathbb{P}(0 \oplus K)$ removed. Moreover, the fiber over $1 \in \Rat(X)$ is the diagonally embedded $\mathbb{P}(K)$. Note that (\ref{ratgrp}) is a monoid homomorphism, so (\ref{ratgrp}) is a bundle with fiber $\mathbb{P}(K)$.

\begin{lemma}
\label{vanish}
We have $\Rat(X)^{\vee} = 1$.
\end{lemma}

\begin{proof}
Suppose we are given a function $(K^{\times} \times K^{\times})/\G_m \times S \to \mathbb{A}^1$. We claim it extends to $\mathbb{P}(K \oplus K) \times S$, which will imply it is pulled back from $S$. Since (\ref{ratgrp}) is a Zariski epimorphism, it will follow that any function on $\Rat(X) \times S$ is pulled back from $S$, and hence $\underline{\Hom}(\Rat(X),\G_m) = 1$. Namely, if $V \subset K$ is a finite-dimensional subspace with $\dim V > 1$, then we obtain a function on the complement in $\mathbb{P}(V \oplus V) \times S$ of $(\mathbb{P}(V \oplus \{ 0 \}) \cup \mathbb{P}(\{ 0 \} \oplus V)) \times S$. Since $\mathbb{P}(V \oplus \{ 0 \})$ and $\mathbb{P}(\{ 0 \} \oplus V))$ are subspaces of codimension at least two, the function extends to $\mathbb{P}(V \oplus V) \times S$ as desired.

Given a multiplicative line bundle on $\Rat(X) \times S$, pull back along the map (\ref{ratgrp}) to obtain a line bundle $\mathcal{L}$ on $(K^{\times} \times K^{\times})/\G_m \times S$. It is enough to show that $\mathcal{L}$ is pulled back from $S$, because (\ref{ratgrp}) is an ind-projective bundle and hence inverse image of line bundles is fully faithful. Lemma \ref{extend} implies that $\mathcal{L}$ extends to $\mathbb{P}(K \oplus K) \times S$. Thus, in order to prove that $\mathcal{L}$ is pulled back from $S$, we need only show that its degree is zero. For this, note that $\mathcal{L}$ is trivial when restricted to the diagonal $\mathbb{P}(K)$, the latter being the kernel of (\ref{ratgrp}).
\end{proof}

\subsection{} Finally, we complete the proof of the main result in the unramified case.

\begin{proof}[Proof of Theorem \ref{mainthm} for $X = \overline{X}$]
The Abel-Jacobi map factors as $X \to \Gr(X) \to \Pic(X)$. By Lemma \ref{globgrdual}, restriction along $X \to \Gr(X)$ induces an isomorphism $\Gr(X)^{\vee} \tilde{\to} \Pic(X)$. Lemma \ref{vanish} then implies that $\Pic(X)^{\vee} \to \Gr(X)^{\vee}$ is an isomorphism.
\end{proof}

In particular, we obtain a perfect pairing \[ \langle \, , \, \rangle_X : \Pic(X) \times \Pic(X) \longrightarrow B\G_m. \]

We also record for later use the following result, which can be interpreted as the unramified case of the local-global compatibility in Theorem \ref{locramthm}.

\begin{proposition}
\label{locglobcompunr}
Let $I$ denote a finite set of closed points in $X$. Then the square \[
\begin{tikzcd}
\Pic(X) \times \prod_{x \in I} \Gr(D_x) \arrow{r} \arrow{d} & \prod_{x \in I} \Pic(D_x) \times \prod_{x \in I} \Gr(D_x) \arrow{d}{\prod_{x \in I} \langle \, , \, ]_{D_x}} \\
\Pic(X) \times \Pic(X) \arrow{r}{\langle \, , \, \rangle_{X}} & B\G_m
\end{tikzcd} \]
commutes up to canonical natural isomorphism.
\end{proposition}

\begin{proof}
This amounts to an isomorphism \[ \langle \mathcal{L},\mathcal{O}_X(Z) \rangle_X \tilde{\longrightarrow} \bigotimes_{x \in I} \langle \mathcal{L}|_{D_x}, Z_x ]_{D_x} \] for any line bundle $\mathcal{L}$ on $X \times S$ and any divisor $Z = \sum_{x \in I} Z_x$ where $Z_x$ is supported in $D_x \times S$. By Theorem \ref{locunrthm} it suffices to treat the case $Z = \Gamma_y$ where $y : S \to D_x$ for some $x \in I$. But then we have \[ \langle \mathcal{L},\mathcal{O}_X(\Gamma_y) \rangle_X \tilde{\longrightarrow} y^*\mathcal{L} \tilde{\longrightarrow} \langle \mathcal{L}|_{D_x}, \Gamma_y ]_{D_x}. \]
\end{proof}

\section{The ramified global Albanese property}

\subsection{} In what follows, let $I$ denote the (finite) set of points in $\overline{X}\setminus X$. We begin with the following general lemma. 

\begin{lemma}
\label{reflexlem}
Let $G$ be a commutative group stack. Let $A$ be a $0$-reflexive sheaf of abelian groups with $\uExt^1(A, \G_m)=\uExt^1(A^*,\G_m)=1$. If there exists an exact triangle
\[ A \longrightarrow H \longrightarrow G \]
where $H$ is $1$-reflexive, then $G$ is $1$-reflexive. 
\end{lemma}

\begin{proof}
Since $A$ is $0$-reflexive and $\uExt^1(A,\G_m) = 1$, we have that $A$ is $1$-reflexive with $A^{\vee} \simeq BA^*$. From the fact that $H^0(BA^*) = 1$, it follows that the dual sequence
\[ G^{\vee} \longrightarrow H^{\vee} \longrightarrow BA^* \]
is exact (Proposition 3.18 in \cite{Br}), hence so is the rotated triangle \[ A^* \longrightarrow G^{\vee} \longrightarrow H^{\vee}. \] Apply the same argument again to obtain the exact triangle \[ A^{**} \longrightarrow H^{\vee \vee} \longrightarrow G^{\vee \vee}. \] Clearly the diagram \[
\begin{tikzcd}
A \arrow{r} \arrow{d} & H \arrow{r} \arrow{d} & G \arrow{d} \\
A^{**} \arrow{r} & H^{\vee \vee} \arrow{r} & G^{\vee \vee}
\end{tikzcd} \]
commutes, so the $1$-reflexivity of $G$ follows from the hypotheses that $A$ is $0$-reflexive and $H$ is $1$-reflexive.
\end{proof}

\subsection{} \label{mainthmproof}
We may now prove Theorem \ref{mainthm} for an open curve. 

\begin{proof}[Proof of Theorem \ref{mainthm}]
First, note that Lemma \ref{reflexlem} applied to the exact triangle
\begin{equation}
\label{restri}
\prod_{x \in I} \Gr(D_x) \longrightarrow \Pic(\overline{X}) \longrightarrow \Pic(X)    
\end{equation}
implies that $\Pic(X)$ is $1$-reflexive. The hypotheses of the lemma are satisfied for $\Gr(D_x)$ by Theorem \ref{locunrthm}.

Consider the square \[
\begin{tikzcd}
\Pic(\overline{X}) \arrow{r} \arrow{d}{\langle \, , \, \rangle_{\overline{X}}} & \prod_{x \in I} \Pic(D_x) \arrow{d}{\prod_{x \in I} \langle \, , \, ]_{D_x}} \\
\Pic(\overline{X})^{\vee} \arrow{r} & \prod_{x \in I} \Gr(D_x)^{\vee},
\end{tikzcd} \]
which commutes by Proposition \ref{locglobcompunr}. Recall that the vertical arrows are isomorphisms by Proposition \ref{locunrthm} and the $X = \overline{X}$ case of Theorem \ref{mainthm}. Passing to the fibers of the horizontal morphisms, we obtain an isomorphism
\begin{equation}
\label{fibiso}
\Pic(X,\partial X) \tilde{\longrightarrow} \Pic(X)^{\vee}.
\end{equation}
In particular, we see that $\Pic(X,\partial X)$ is also $1$-reflexive.

It remains only to show that the isomorphism \[ \Pic(X,\partial X)^{\vee} \tilde{\longrightarrow} \Pic(X) \] dual to (\ref{fibiso}) agrees with restriction along $\AJ_X$. This amounts to commutativity of the triangle \[
\begin{tikzcd}
X \arrow{r}{\AJ_X} \arrow{dr} & \Pic(X,\partial X) \arrow{d}{(\ref{fibiso})} \\
& \Pic(X)^{\vee},
\end{tikzcd} \]
where the diagonal arrow is the tautological map $y \mapsto (\mathcal{L} \mapsto y^*\mathcal{L})$. Denoting by \[ \langle \, , \, ]_X : \Pic(X) \times \Pic(X,\partial X) \longrightarrow B\G_m \] the pairing corresponding to (\ref{fibiso}), we must exhibit a natural isomorphism \[ \langle (\mathcal{O}_{\overline{X}}(y),1),\mathcal{L} ]_X \tilde{\longrightarrow} y^*\mathcal{L} \] for any $y : S \to X$ and any extendable line bundle $\mathcal{L}$ on $X \times S$.

First, for $\mathcal{L} = \overline{\mathcal{L}}|_X$ where $\overline{\mathcal{L}}$ is a line bundle on $\overline{X} \times S$, we have \[ \langle (\mathcal{O}_{\overline{X}}(y),1),\overline{\mathcal{L}}|_X ]_X = \langle \mathcal{O}_{\overline{X}}(y),\overline{\mathcal{L}} \rangle_{\overline{X}} \tilde{\longrightarrow} y^*\overline{\mathcal{L}} = y^*\mathcal{L}. \] In view of the exact triangle (\ref{restri}), it remains to show that in the case $\overline{\mathcal{L}} = \mathcal{O}_{\overline{X}}(Z)$ where $Z = \sum_{x \in I} Z_x$ with $Z_x$ supported on $\sqcup_{x \in I} D_x \times S$, the resulting trivializations of both sides of \[ \langle \mathcal{O}_{\overline{X}}(y),\overline{\mathcal{L}} \rangle_{\overline{X}} \tilde{\longrightarrow} y^*\mathcal{L} \] match up. Here the left side is trivialized via \[ \langle \mathcal{O}_{\overline{X}}(y),\mathcal{O}_{\overline{X}}(Z) \rangle_{\overline{X}} \tilde{\longrightarrow} \bigotimes_{x \in I} \langle \mathcal{O}_{\overline{X}}(y)|_{D_x},Z_x ]_{D_x} \tilde{\longrightarrow} \bigotimes_{x \in I} \langle \mathcal{O}_{D_x},Z_x ]_{D_x} \tilde{\longrightarrow} \mathcal{O}_S. \] One verifies by a diagram chase that this coincides with the trivialization obtained by restricting $1 : \mathcal{O}_{\overline{X}} \to \mathcal{O}_{\overline{X}}(Z)$ along $y$.

\end{proof}

\begin{proof}[Proof of Corollary \ref{maincor}]
The proof is identical to that of Corollary \ref{locramalb}, replacing $\mathfrak{L}_x\G_m$ by $\Pic(X,\partial X)$ and $\Pic(\mathring{D}_x)$ by $\Pic(X)$.
\end{proof}

\subsection{} \label{locramthmproof}
We now complete the proof of Theorem \ref{locramthm} by establishing local-global compatibility.

\begin{proof}[Proof of Theorem \ref{locramthm}]
The nondegeneracy of $\langle \, , \, ]_{\mathring{D}_{\Ran}}$ is Theorem \ref{locnondeg}, and compatibility with factorization structures is clear from the construction of the pairing. Thus it remains only to prove the commutativity of the diagram (\ref{locglobdiag}). Commutativity of the right square amounts to commutativity of the square \[
\begin{tikzcd}
\Pic(X,\partial X) \times \Pic(X,\partial X) \arrow{r} \arrow{d} & \Pic(X) \times \Pic(X,\partial X) \arrow{d}{\langle \, , \, ]_X} \\
\Pic(X,\partial X) \times \Pic(X)\arrow{r}{[ \, , \, \rangle_X} & B\G_m.
\end{tikzcd} \]
This is proved similarly to Proposition \ref{locpairingcomm}: by Theorem \ref{mainthm} it suffices to show the two compositions agree when restricted along $\AJ_X \times \AJ_X$, and indeed both restrictions identify with $\mathcal{O}_{X \times X}(\Delta)$.

Finally, we establish commutativity of the left square in (\ref{locglobdiag}), or equivalently commutativity of the square  \[
\begin{tikzcd}
\Pic(X) \times \prod_{x \in I} \mathfrak{L}_x\G_m \arrow{r} \arrow{d} & \prod_{x \in I} \Pic(\mathring{D}_x) \times \prod_{x \in I} \mathfrak{L}_x\G_m \arrow{d}{\prod_{x \in I} \langle \, , \, ]_{\mathring{D}_x}} \\
\Pic(X) \times \Pic(X,\partial X) \arrow{r}{\langle \, , \, ]_X} & B\G_m.
\end{tikzcd} \]
Recall that we have an exact triangle \[ \prod_{x \in I} \Gr(D_x) \longrightarrow \Pic(\overline{X}) \longrightarrow \Pic(X), \] and hence an exact triangle \[ \Pic(X)^{\vee} \longrightarrow \Pic(\overline{X})^{\vee} \longrightarrow \prod_{x \in I} \Gr(D_x)^{\vee} \] because $\underline{\Ext}^1(\Gr(D_x),\G_m) = 1$ by Theorem \ref{locunrthm}. Thus it suffices to construct an isomorphism of line bundles
\begin{equation}
\label{compiso}
\langle \mathcal{L}|_X,(\mathcal{O}_X,f) ]_X \tilde{\longrightarrow} \bigotimes_{x \in I} \langle \mathcal{L}|_{\mathring{D}_x},f_x ]_{\mathring{D}_x}
\end{equation} for any line bundle $\mathcal{L}$ on $\overline{X} \times S$ and any $f = (f_x)_{x \in I}$ where $f_x : \mathring{D}_{x,S} \to \G_m$, such that for $\mathcal{L} = \mathcal{O}_{\overline{X}}(Z)$ with $Z$ a divisor supported on $\sqcup_{x \in I} D_x \times S$, the resulting trivializations of both sides match up under (\ref{compiso}).

Indeed, we have
\begin{align*}
\langle \mathcal{L}|_X,(\mathcal{O}_X,f) ]_X &\tilde{\longrightarrow} \langle \mathcal{L},\mathcal{O}_{\overline{X}}(\ddiv(f)) \rangle_{\overline{X}} \\
&\tilde{\longrightarrow} \bigotimes_{x \in I} \langle \mathcal{L}|_{D_x},\ddiv(f_x) ]_{D_x} \\
&\tilde{\longrightarrow} \bigotimes_{x \in I} \langle \mathcal{L}|_{\mathring{D}_x},f_x ]_{\mathring{D}_x}.
\end{align*}
Here the first isomorphism was established in the proof of Theorem \ref{mainthm}, the second isomorphism is Proposition \ref{locglobcompunr}, and the last isomorphism is by the definition of $\langle \, , \, ]_{D_x}$.

For $\mathcal{L} = \mathcal{O}_{\overline{X}}(Z)$, one traces through the constructions to see that the resulting trivializations \[ \langle \mathcal{L}|_X,(\mathcal{O}_X,f) ]_X \tilde{\longrightarrow} \langle \mathcal{O}_X,(\mathcal{O}_X,f) ]_X \tilde{\longrightarrow} \mathcal{O}_S \] and \[ \langle \mathcal{L}|_{\mathring{D}_x},f_x ]_{\mathring{D}_x} \tilde{\longrightarrow} \langle \mathcal{O}_{\mathring{D}_x},f_x ]_{\mathring{D}_x} \tilde{\longrightarrow} \mathcal{O}_S \] match up under (\ref{compiso}).
\end{proof}

\section{Global determinant formulas}

\subsection{} The goal of this section is to give a formula for the global pairing \[ \langle \, , \, ]_X : \Pic(X) \times \Pic(X,\partial X) \longrightarrow B\G_m \] in terms of determinants, analogously to how the local pairing $\langle \, , \, ]_{\mathring{D}_{\Ran}}$ was defined.

First we treat the case where $X = \overline{X}$ is projective. We will need the following lemma. Fix $x_I : S \to X^I$ and a line bundle $\mathcal{L}$ on $X \times S$.

\begin{lemma}
\label{colatlem}
Restriction from $X \times S$ to $\mathring{D}_{x_I}$ realizes $\Gamma((X \times S) \setminus \Gamma_{x_I},\mathcal{L})$ as a colattice in the Tate vector bundle $\Gamma(\mathring{D}_{x_I},\mathcal{L})$ on $S$.
\end{lemma}

\begin{proof}
Twisting by a divisor supported on $\Gamma_{x_I}$, we may assume that $\mathcal{L}$ is of arbitrarily large degree, and hence that $R^1\Gamma(X \times S,\mathcal{L}) = 0$. Thus $R^0\Gamma(X \times S,\mathcal{L})$ is a vector bundle on $S$, and we have a short exact sequence \[ 0 \longrightarrow R^0\Gamma(X \times S,\mathcal{L}) \longrightarrow \Gamma((X \times S) \setminus \Gamma_{x_I},\mathcal{L}) \oplus \Gamma(D_{x_I},\mathcal{L}) \longrightarrow \Gamma(\mathring{D}_{x_I},\mathcal{L}) \longrightarrow 0. \] Since $\Gamma(D_{x_I},\mathcal{L})$ is linearly compact and open in $\Gamma(\mathring{D}_{x_I},\mathcal{L})$, the claim follows.
\end{proof}

Continuing to assume that $X = \overline{X}$, for any line bundle $\mathcal{L}$ on $X \times S$ we obtain a perfect complex $R\Gamma(X \times S,\mathcal{L})$ on $S$, which has a well-defined determinant line.

\begin{theorem}
\label{globaldetformula}
For any line bundles $\mathcal{L}$ and $\mathcal{M}$ on $X \times S$, the line bundle $\langle \mathcal{L},\mathcal{M} \rangle_X$ on $S$ is naturally isomorphic to \[ \det R\Gamma(X \times S,\mathcal{L} \otimes \mathcal{M}) \otimes (\det R\Gamma(X \times S,\mathcal{L}))^{-1} \otimes (\det R\Gamma(X \times S,\mathcal{M}))^{-1} \otimes \det R\Gamma(X \times S,\mathcal{O}). \]
\end{theorem}

\begin{proof}
We are comparing two line bundles on the product $\Pic(X) \times \Pic(X)$. It was shown in the proof of Lemma \ref{vanish} that any function on $\Rat(X)$ is constant, so it suffices to identify these line bundles after restriction along \[ \Pic(X) \times \Gr(X) \longrightarrow \Pic(X) \times \Pic(X). \] Furthermore, we saw in the proof of Lemma \ref{globgrdual} that \[ |\Ran^{\bullet} \times \Gr(D)_{\Ran}| \tilde{\longrightarrow} \Gr(X), \] whence it is enough to produce a $\Ran$-equivariant isomorphism between the two line bundles after restriction along \[ \Pic(X) \times \Gr(D)_{\Ran} \longrightarrow \Pic(X) \times \Gr(X). \]

Fix $x_I : S \to X^I$, a line bundle $\mathcal{L}$ on $X \times S$, and a divisor $Z$ relative to $S$ supported on $D_{x_I}$. Recall that we have a canonical isomorphism \[ \langle \mathcal{L},\mathcal{O}(Z) \rangle_X \tilde{\longrightarrow} \langle \mathcal{L}|_{D_{x_I}},Z ]_{D_{\Ran}} \] (cf. Proposition \ref{locglobcompunr}).

In order to relate $\langle \mathcal{L}|_{D_{x_I}},Z ]_{D_{\Ran}}$ to the desired determinant formula, we use a \v{C}ech resolution in the form of the triangle \[ R\Gamma(X \times S,\mathcal{L}) \longrightarrow \Gamma((X \times S) \setminus \Gamma_{x_I},\mathcal{L}) \oplus \Gamma(D_{x_I},\mathcal{L}) \longrightarrow \Gamma(\mathring{D}_{x_I},\mathcal{L}), \] which is exact by the Beauville-Laszlo theorem. By Lemma \ref{colatlem} we obtain an isomorphism
\begin{equation}
\label{gluedet}
\det(\Gamma((X \times S) \setminus \Gamma_{x_I},\mathcal{L}),\Gamma(D_{x_I},\mathcal{L})) \tilde{\longrightarrow} \det R\Gamma(X \times S,\mathcal{L}).
\end{equation}

Recall that by Proposition \ref{unrpairing}, we have \[ \langle \mathcal{L}|_{D_{x_I}},Z]_{D_{\Ran}} = \det(\Gamma(D_{x_I},\mathcal{L}),\Gamma(D_{x_I},\mathcal{L}(Z)))^{-1} \otimes \det(\Gamma(D_{x_I},\mathcal{O}),\Gamma(D_{x_I},\mathcal{O}(Z))). \] Using the aforementioned \v{C}ech resolution, we find that
\begin{align*}
\det(\Gamma(D_{x_I},&\mathcal{L}),\Gamma(D_{x_I},\mathcal{L}(Z))) \\ &= \det(\Gamma(D_{x_I},\mathcal{L}),\Gamma((X \times S) \setminus \Gamma_{x_I},\mathcal{L})) \otimes \det(\Gamma((X \times S) \setminus \Gamma_{x_I},\mathcal{L}),\Gamma(D_{x_I},\mathcal{L}(Z))) \\
&= \det(\Gamma(D_{x_I},\mathcal{L}),\Gamma((X \times S) \setminus \Gamma_{x_I},\mathcal{L})) \otimes \det(\Gamma((X \times S) \setminus \Gamma_{x_I},\mathcal{L}(Z)),\Gamma(D_{x_I},\mathcal{L}(Z))) \\
&= \det(R\Gamma(X \times S,\mathcal{L})) \otimes (\det R\Gamma(X \times S,\mathcal{L}(Z)))^{-1}.
\end{align*}

Thus we have given the desired formula for the pairing on $\Pic(X) \times \Gr(D)_{\Ran}$. The $\Ran$-equivariance follows from the compatibility of the \v{C}ech resolution with refinement of the covering.
\end{proof}

\subsection{} We now introduce a formula for $\langle \, , \, ]_X$ in the case where $X$ is affine. Let $I$ denote the set of closed points in $\overline{X} \setminus X$. We will use the notations \[ D_{I,S} := \bigsqcup_{x \in I} D_{x,S} \text{ and } \mathring{D}_{I,S} := \bigsqcup_{x \in I} \mathring{D}_{x,S}. \]

Fix $S$-points of $\Pic(X)$ and $\Pic(X,\partial X)$, which are respectively an extendable line bundle $\mathcal{L}$ on $X \times S$ and a line bundle $\mathcal{M}$ on $X \times S$ together with a trivialization $\alpha$ of $\mathcal{M}$ over $\mathring{D}_{I,S}$. We obtain an isomorphism of Tate vector bundles \[ \Gamma(\mathring{D}_{I,S},\mathcal{L}) \xrightarrow{\id_{\mathcal{L}} \otimes \alpha} \Gamma(\mathring{D}_{I,S},\mathcal{L} \otimes \mathcal{M}). \] In particular, by comparing under this isomorphism we may form the relative determinant of the colattices $\Gamma(X \times S,\mathcal{L})$ and $\Gamma(X \times S,\mathcal{L} \otimes \mathcal{M})$.

\begin{theorem}
For $X$ affine, the line bundle $\langle \mathcal{L},(\mathcal{M},\alpha) ]_X$ on $S$ is naturally isomorphic to \[ \det\big((\id_{\mathcal{L}} \otimes \alpha)(\Gamma(X \times S,\mathcal{L})),\Gamma(X \times S,\mathcal{L} \otimes \mathcal{M})\big)^{-1} \otimes \det\big(\alpha(\Gamma(X \times S,\mathcal{O})),\Gamma(X \times S,\mathcal{M})\big). \]
\end{theorem}

\begin{proof}
Since we are trying to compare two line bundles on $\Pic(X,\partial X) \times \Pic(X)$, it therefore suffices to construct a $\Gr(D_I)$-equivariant isomorphism after restriction along the $\Gr(D_I)$-torsor \[ \Pic(X,\partial X) \times \Pic(\overline{X}) \longrightarrow \Pic(X,\partial X) \times \Pic(X). \]

So suppose we are given an extension $\overline{\mathcal{L}}$ of $\mathcal{L}$ to $\overline{X} \times S$. We also write $\overline{\mathcal{M}}$ for the the extension of $\mathcal{M}$ obtained using $\alpha$, so that $\alpha$ extends to a trivialization of $\overline{\mathcal{M}}$ over $D_{I,S}$. The isomorphism \ref{gluedet} therefore gives rise to an isomorphism
\begin{align*}
\det\big((\id_{\mathcal{L}} \otimes \alpha)(\Gamma(&X \times S,\mathcal{L})),\Gamma(X \times S,\mathcal{L} \otimes \mathcal{M})\big) \\
&\tilde{\longrightarrow} \det(R\Gamma(\overline{X} \times S,\overline{\mathcal{L}})) \otimes \det(R\Gamma(\overline{X} \times S,\overline{\mathcal{L}} \otimes \overline{\mathcal{M}}))^{-1}.
\end{align*}

Recall that in the proof of Theorem \ref{mainthm}, we established the commutativity of the square \[
\begin{tikzcd}
\Pic(X,\partial X) \times \Pic(\overline{X}) \arrow{r} \arrow{d} & \Pic(X,\partial X) \times \Pic(X) \arrow{d}{\langle \ , \ ]_X} \\
\Pic(\overline{X})\times \Pic(\overline{X}) \arrow{r}{\langle \ , \ \rangle_{\overline{X}}} & B\G_m.
\end{tikzcd} \]
Thus, by Theorem \ref{globaldetformula} we obtain the desired isomorphism of line bundles on $\Pic(X,\partial X) \times \Pic(X)$. Tracing through the construction of the isomorphism, one finds that it is $\Gr(D_I)$-equivariant as needed.
\end{proof}

\section{Geometric class field theory}

\label{cftsec}

\subsection{} In this section we prove Theorem \ref{finext}. We will need a couple of technical lemmas. Let $p$ denote the characteristic of $k$. Fix a closed point $x$ in $X$.

\begin{lemma}
\label{tamelemma}
For $n$ prime to $p$, the pullback map 
\[ H^1(S,\G_m)[n] \longrightarrow H^1(\mathring{D}_{x,S},\G_m)[n] \]
is an isomorphism.
\end{lemma}

\begin{proof}
The multiplication by $n$ map $\Z \to \Z$ induces a long exact sequence on fppf (here equivalent to \'{e}tale) cohomology
\[ 0 \longrightarrow H^0(S,\Z) \longrightarrow H^0(S,\Z) \longrightarrow H^0(S,\Z/n \Z) \longrightarrow H^1(S,\Z) \longrightarrow \cdots \]
Since for $A$ a discrete abelian group $H^0(S,A)$ classifies locally constant maps $S \to A$, the map 
\[ H^0(S,\Z) \longrightarrow H^0(S,\Z/n \Z) \]
is surjective. Hence $H^1(S,\Z)$ is $n$-torsion free. 

Now by Section $3.1$ in \cite{BC}, we have an isomorphism
\[ H^1(\mathring{D}_{x,S},\G_m) \tilde{\longrightarrow} H^1(S,\Z) \oplus H^1(\mathbb{A}^1 \times S,\G_m)_0 \oplus H^1(S,\G_m), \]
where $H^1(\mathbb{A}^1 \times S,\G_m)_0$ is the subgroup of line bundles on $\mathbb{A}^1 \times S$ whose restriction to $\{0\} \times S$ is trivializable. We have seen that $H^1(S,\Z)[n] = 0$, so it remains to show that $H^1(\mathbb{A}^1 \times S,\G_m)_0$ is $n$-torsion free. The boundary map attached to the Kummer sequence yields a surjection \[ H^1(\A^1 \times S,\mu_n) \longrightarrow H^1(\A^1 \times S,\G_m)[n], \] and since $p$ does not divide $n$, we have $\mu_n \cong \Z/n\Z$ and \[ H^1(S,\mu_n) \tilde{\longrightarrow} H^1(\A^1 \times S,\mu_n). \] It follows that $H^1(\A^1 \times S,\G_m)_0[n] = 1$.
\end{proof}

\begin{lemma}
\label{wildlemma}
The inclusion
\[ \uHom(\mu_{p^k}, \Pic(\mathring{D}_x)) \longrightarrow \uHom(\mu_{p^k},\uMap(\mathring{D}_x,B\G_m)) \]
is an isomorphism.
\end{lemma}

\begin{proof}

A homomorphism $\mu_{p^k} \times S \to \uMap(\mathring{D}_x,B\G_m) \times S$ over $S$ is in particular a map $\mu_{p^k} \times S \to \uMap(\mathring{D}_x,B\G_m)$ whose restriction to $\{1\} \times S$ is the trivial bundle. This amounts to a line bundle on $\mathring{D}_{x,\mu_{p^k} \times S}$ equipped with a trivialization over $\mathring{D}_{x,\{ 1 \} \times S}$. Since \[ \mathring{D}_{x,S} \longrightarrow \mathring{D}_{x,\mu_{p^k} \times S} \] is a nilpotent embedding and $\G_m$ is smooth, the trivialization extends. Hence the bundle is trivializable, and in particular extendable. 
\end{proof}

\subsection{} \label{finextproof}
We are now in a position to prove the theorem.

\begin{proof}[Proof of Theorem \ref{finext}]
Note that the isomorphism for $\mathring{D}_x$ implies the one for $X$. Letting $I$ denote the set of closed points in $\overline{X} \setminus X$, this follows from the Cartesian square \[
\begin{tikzcd}
\Pic(X) \arrow{r} \arrow{d} & \uMap(X,B\G_m) \arrow{d} \\
\prod_{x \in I} \Pic(\mathring{D}_x) \arrow{r} & \prod_{x \in I} \uMap(\mathring{D}_x,B\G_m).
\end{tikzcd} \]

One also reduces immediately to the case that $A = \Z/n\Z$ is cyclic. Thus we are trying to prove that \[ \uHom(\mu_n,\Pic(\mathring{D}_x)) \longrightarrow \uHom(\mu_n,\Map(\mathring{D}_x,B\G_m)) \] is an isomorphism. If $n$ is prime to $p$, then $\mu_n \cong \Z/n\Z$ and the claim follows from Lemma \ref{tamelemma}. Then it remains only to treat the case $n = p^k$, which is precisely Lemma \ref{wildlemma}.
\end{proof}

\appendix
\section{The ramified local Abel-Jacobi map}

\label{appendix}

\subsection{} Let us fix a closed point $x$ in $X$ once and for all. For simplicity we assume that $x$ is a $k$-rational point. We will suppress it from the notation and write $D := D_x$, $\mathring{D} := \mathring{D}_x$, etc.

We now recall the construction of the ramified local Abel-Jacobi map introduced by Contou-Carr\`{e}re in \cite{CC1}. Viewing $\mathfrak{L}\G_m$ as a trivializable $\mathfrak{L}^+\G_m$-torsor over $\Gr$, we will now construct a \emph{canonical} extension of $\mathfrak{L}\G_m$ to a trivializable $\mathfrak{L}^+\G_m$-torsor $\mathfrak{L}\G_m^{\aff}$ over $\aff(\Gr)$ (recall that $\aff(Y) := \Gamma(Y,\mathcal{O}_Y)$ for $Y$ a prestack).

Let $\mathcal{L}_{\univ}$ denote the tautological line bundle on $D \times \Gr$. The Weil restriction of $\mathcal{L}_{\univ}$ along the projection $D \times \Gr \to \Gr$ is precisely the $\mathfrak{L}^+\G_m$-torsor $\mathfrak{L}\G_m$. Observe that $\mathcal{L}_{\univ}$ canonically extends to a line bundle $\widetilde{\mathcal{L}}_{\univ}$ on $\aff(D \times \Gr)$: since $\mathcal{L}_{\univ}$ is trivializable, the $\Gamma(D \times \Gr,\mathcal{O})$-module $\Gamma(D \times \Gr,\mathcal{L}_{\univ})$ is free of rank one and we let $\widetilde{\mathcal{L}}_{\univ}$ be the associated line bundle on $\aff(D \times \Gr)$. Let $\mathcal{L}_{\univ}^{\aff}$ be the inverse image of $\widetilde{\mathcal{L}}_{\univ}$ along the map \[ D \times \aff(\Gr) \longrightarrow \aff(D \times \aff(\Gr)) \tilde{\longrightarrow} \aff(D \times \Gr). \] Finally, we define $\mathfrak{L}\G_m^{\aff}$ as the Weil restriction of $\mathcal{L}_{\univ}^{\aff}$ along the projection $D \times \aff(\Gr) \to \aff(\Gr)$.

\subsection{} Observe that the inverse image of $\mathcal{L}_{\univ}$ along the map \[ \id_D \times \AJ_D : D \times D \to D \times \Gr \] canonically identifies with $\mathcal{O}_{D \times D}(\Delta)$. In particular, we have a canonical section \[ \sigma_{\can} \in \Gamma(D \times D,\mathcal{L}_{\univ}) = \Gamma(D \times \aff(D),\mathcal{L}_{\univ}^{\aff}), \] where the equality here follows from the observation that $\aff(D \times \aff(D)) = \aff(D \times D)$.

\begin{proposition}

The section $\sigma_{\can}$ of $\mathcal{L}_{\univ}^{\aff}$ is nonvanishing over $D \times \mathring{D}$.

\end{proposition}

\begin{proof}

The section $\sigma_{\can}$ corresponds to the canonical morphism $\mathcal{O}_{D \times D} \to \mathcal{O}_{D \times D}(\Delta)$, which is dual to the inclusion $I_{\Delta} \to \mathcal{O}_{D \times D}$ of the ideal of the diagonal. Choosing a coordinate $t$ on $D$, the claim now follows from the observation that $t-u$ is a unit in \[ \Gamma(D \times \mathring{D},\mathcal{O}_{D \times \mathring{D}}) \cong k((t))[[u]] \] because
\begin{equation}
\label{geomser}
\frac{1}{t-u} = \frac{1}{t} + \frac{u}{t^2} + \frac{u^2}{t^3} + \cdots.
\end{equation}

\end{proof}

Thus $\sigma_{\can}$ determines a trivialization of $\mathcal{L}_{\univ}^{\aff}$ over $D \times \mathring{D}$, and hence a section of its Weil restriction $\mathfrak{L}\G_m^{\aff}$ over $\mathring{D}$. The \emph{ramified local Abel-Jacobi map} is the resulting map \[ \AJ_{\mathring{D}} : \mathring{D} \longrightarrow \mathfrak{L}\G_m^{\aff}. \] The above constructions can be performed relative to an affine base scheme $S$: we have the line bundle $(\mathcal{L}_{\univ} \boxtimes \mathcal{O}_S)^{\aff}$ on $D \times \aff(\Gr \times S)$, as well as its Weil restriction $(\mathfrak{L}\G_m \times S)^{\aff}$, an $\mathfrak{L}^+\G_m$-bundle on $\aff(\Gr \times S)$. We denote by \[ \AJ_{\mathring{D}_S} : \mathring{D}_S \longrightarrow (\mathfrak{L}\G_m \times S)^{\aff} \] the resulting parameterized version of the local Abel-Jacobi map.

By construction, we have a commutative square
\begin{equation}
\label{lajsq}
\begin{tikzcd}
\mathring{D}_S \arrow{r}{\AJ_{\mathring{D}_S}} \arrow{d} & (\mathfrak{L}\G_m \times S)^{\aff} \arrow{d} \\
\aff(D \times S) \arrow{r}[yshift=1ex]{\aff(\AJ_D \times \id_S)} & \aff(\Gr \times S).
\end{tikzcd}
\end{equation}

\subsection{} The operation of ``restriction along the local Abel-Jacobi map" requires some care to define. First, observe that for any affine scheme $Y$, the restriction map \[ \Map((\mathfrak{L}\G_m \times S)^{\aff},Y) \longrightarrow \Map(\mathfrak{L}\G_m \times S,Y) \] is injective. This reduces immediately to the case $Y = \A^1$, where it follows from the observation that functions on $\mathfrak{L}\G_m \times S$ are the completion of functions on $(\mathfrak{L}\G_m \times S)^{\aff}$ with respect to a separated filtration.

\begin{proposition}
If $G$ is a commutative affine group scheme, then the image of the map \[ \Hom_S(\mathfrak{L}\G_m \times S,G \times S) \longrightarrow \Map(\mathfrak{L}\G_m \times S,G) \] is contained in the image of $\Map((\mathfrak{L}\G_m \times S)^{\aff},G)$.
\end{proposition}
\begin{proof}
Choosing a splitting $\mathfrak{L}\G_m \cong \mathfrak{L}^+\G_m \times \Gr$, one checks that the composition
\begin{align*}
\Hom_S(\mathfrak{L}\G_m \times S,G \times S) &\tilde{\longrightarrow} \Hom_S(\mathfrak{L}^+\G_m \times S,G \times S) \times \Hom_S(\Gr \times S,G \times S) \\
&\longrightarrow \Map(\mathfrak{L}^+\G_m \times S,G) \times \Map(\Gr \times S,G) \\
&\tilde{\longrightarrow} \Map(\mathfrak{L}^+\G_m \times S,G) \times \Map(\aff(\Gr \times S),G) \\
&\longrightarrow \Map((\mathfrak{L}\G_m \times S)^{\aff},G) \\
&\longrightarrow \Map(\mathfrak{L}\G_m \times S,G)
\end{align*}
is the map in question. Here the second to last map in the composition sends $f : \mathfrak{L}^+\G_m \times S \to G$ and $g : \aff(\Gr \times S) \to G$ to
\begin{align*}
(\mathfrak{L}\G_m \times S)^{\aff} &\tilde{\longrightarrow} (\mathfrak{L}^+\G_m \times S) \times_S \aff(\Gr \times S) \\\
&\longrightarrow \mathfrak{L}^+\G_m \times S \times \aff(\Gr \times S) \\
&\xrightarrow{f \times g} G \times G \longrightarrow G.
\end{align*}
\end{proof}

Thus we obtain a canonical map
\begin{equation}
\label{oblvalb}
\Hom_S(\mathfrak{L}\G_m \times S,G \times S) \longrightarrow \Map((\mathfrak{L}\G_m \times S)^{\aff},G).
\end{equation}

Finally, we consider the composition \[ \Hom_S(\mathfrak{L}\G_m \times S,G \times S) \xrightarrow{(\ref{oblvalb})} \Map((\mathfrak{L}\G_m \times S)^{\aff},G) \longrightarrow \Map(\mathring{D}_S,G), \] where the second map is precomposition with $\AJ_{\mathring{D}_S}$. This is easily seen to be natural in $S$ and hence induces
\begin{equation}
\label{localb}
\underline{\Hom}(\mathfrak{L}\G_m,G) \longrightarrow \underline{\Map}(\mathring{D},G).
\end{equation}

On the other hand, we have the perfect pairing (\ref{ccpairing}) arising as the commutator in the Tate extension.

\begin{proposition}
When $G = \G_m$, the map (\ref{localb}) is inverse to the isomorphism induced by the pairing (\ref{ccpairing}).
\end{proposition}

\begin{proof}
First observe that we have an isomorphism of short exact sequences \[
\begin{tikzcd}
\mathfrak{L}^+\G_m \arrow{r} \arrow{d} & \mathfrak{L}\G_m \arrow{r} \arrow{d}{(\ref{ccpairing})} & \Gr \arrow{d} \\
\Gr^* \arrow{r} & \mathfrak{L}\G_m^* \arrow{r}  & \mathfrak{L}^+\G_m^*,
\end{tikzcd} \]
where the left and right vertical isomorphisms are induced by the pairing $[ \, , \, \rangle_D$. Namely, this results from the diagram in the proof of Theorem \ref{locnondeg} by restricting to $x$ and rotating the rows.

To prove the proposition, it therefore suffices to prove that the diagram \[
\begin{tikzcd}
\Gr^* \arrow{r} \arrow{d} & \mathfrak{L}\G_m^* \arrow{r} \arrow{d}{(\ref{localb})}  & \mathfrak{L}^+\G_m^* \arrow{d} \\
\mathfrak{L}^+\G_m \arrow{r} & \mathfrak{L}\G_m \arrow{r} & \Gr
\end{tikzcd} \]
commutes, where the left and right vertical maps are inverse to those in the previous diagram. By Theorem \ref{locunrthm}, the left vertical map agrees with restriction along $\AJ_D$, and therefore the left square commutes because the square (\ref{lajsq}) does.

As for the right square, fix a coordinate $t$ at $x$, which induces a splitting of the projection $\mathfrak{L}\G_m \to \Gr$ and hence $\mathfrak{L}\G_m^{\aff} \cong \aff(\Gr) \times \mathfrak{L}^+\G_m$. Observe that the resulting composition
\begin{equation}
\label{lajproj}
\mathring{D} \xrightarrow{\AJ_{\mathring{D}}} \mathfrak{L}\G_m^{\aff} \longrightarrow \mathfrak{L}^+\G_m
\end{equation} is the invertible function on $D \times \mathring{D}$ given by the series $t-u$ (as we previously observed, this is a unit with inverse given by the formula (\ref{geomser})). The claim now amounts to showing that the composition 
\begin{equation}
\label{lajcomp}
\Gr \xrightarrow{[ \, , \, \rangle_D} \mathfrak{L}^+\G_m^* \longrightarrow \mathfrak{L}\G_m \longrightarrow \Gr
\end{equation}
is the identity, where the second map is given by restriction along (\ref{lajproj}). Note that Theorem \ref{locunrthm} implies that restriction along $\AJ_D$ induces an isomorphism \[ \uHom(\Gr,\Gr) \tilde{\longrightarrow} \uMap(D,\Gr), \] so it suffices to check that (\ref{lajcomp}) is the identity after precomposition with $\AJ_D$. Given $f : S \to D$, which is just a nilpotent function on $S$, the resulting element of $\mathfrak{L}^+\G_m^*$ is given by evaluation at $f$. Precomposing with (\ref{lajcomp}) yields the function $t - f : \mathring{D}_S \to \G_m$, whose image in $\Gr$ is $\AJ_D(f)$ as desired.
\end{proof}

\end{document}